\documentclass[12pt]{article}
\usepackage{amsmath,amssymb}
\usepackage{amsthm}
\usepackage{graphicx}
\usepackage{subfigure}
\usepackage{url}
\usepackage{multirow}
\usepackage[font=footnotesize]{subfig}
\usepackage{array}
\usepackage{mathdots}
\usepackage{cellspace}
\usepackage[round]{natbib}
\newtheorem{theorem2}{Theorem}
\newtheorem{proposition}{Proposition}

\begin{document}

\newcommand{\x}{\bx}
\newcommand{\f}{\bdf}
\newcommand{\e}{\mathbf{e}}
\newcommand{\W}{\mathbf{W}}
\newcommand{\R}{\mathbb{R}}
\newcommand{\E}{\mathbb{E}}
\newcommand{\V}{\mathbf{V}}
\newcommand{\I}{\mathbf{I}}
\newcommand{\Z}{\mathbf{Z}}
\newcommand{\X}{\mathbf{X}}
\newcommand{\Y}{\mathbf{Y}}
\newcommand{\T}{\mathbf{T}}
\newcommand{\U}{\mathbf{U}}
\newcommand{\D}{\mathbf{D}}
\newcommand{\dif}{\mathrm{d}}
\renewcommand{\P}{\mathbb{P}}

\newcommand\bx{\mathbf{x}}
\newcommand\bz{\mathbf{z}}
\newcommand\bSi{\mathbf{\Sigma}}
\newcommand\bmu{\boldsymbol{\mu}}
\newcommand\bGa{\mathbf{\Gamma}}
\newcommand\bPsi{\mathbf{\Psi}}
\newcommand\bLa{\mathbf{\Lambda}}
\newcommand\bdf{\mathbf{f}}
\newcommand\bde{\mathbf{e}}
\newcommand\bI{\mathbf{I}}
\newcommand\bS{\mathbf{S}}
\newcommand\bun{\mathbf{1}}
\newcommand\bzero{\mathbf{0}}
\newcommand\cov{\mathop{\text{cov}}}
\newcommand\si{\sigma}

\begin{center}
{\bf \large On estimation of the noise variance in high-dimensional probabilistic principal component analysis}
\\[3mm] Damien {\sc Passemier} \\  Department of Electronic and Computer Engineering \\Hong Kong University of Science and Technology
\\[3mm] Zhaoyuan {\sc Li} \\Department of Statistics and Actuarial Science\\The University of Hong Kong
\\[3mm] Jianfeng {\sc Yao} \\Department of Statistics and Actuarial Science\\The University of Hong Kong
\end{center}

\begin{abstract}
In this paper, we develop new statistical theory for probabilistic principal component analysis models in high dimensions. The focus is the estimation of the noise variance, which is an important and unresolved issue when the number of variables is large in comparison with the sample size. We first unveil the reasons of a widely observed downward bias of the maximum likelihood estimator of the variance when the data dimension is high. We then propose a bias-corrected estimator using random matrix theory and establish its asymptotic normality. The superiority of the new (bias-corrected) estimator over existing alternatives is first checked by Monte-Carlo experiments with various combinations of $(p, n)$ (dimension and sample size). In order to demonstrate further potential benefits from the results of the paper to general probability PCA analysis, we provide evidence of net improvements in two popular procedures \citep{Solo, Baij} for determining the number of principal components when the respective variance estimator proposed by these authors is replaced by the bias-corrected estimator. The new estimator is also used to derive new asymptotics for the related goodness-of-fit statistic under the high-dimensional scheme.

\end{abstract}

\section{Introduction}
Principal component analysis (PCA) is a very popular technique in multivariate analysis for dimensionality reduction and feature extraction. Due to dramatic development in data-collection technology, high-dimensional data are nowadays common in many fields. Natural high-dimensional data, such as images, signal processing, documents and biological data often reside in a low-dimensional subspace or low-dimensional manifold \citep{ding2011bayesian}. In financial econometrics, it is commonly believed that the variations in a large number of economic variables can be modeled by a small number of reference variables \citep{stock1989new, forni2000reference, Baij, BaiJ03}. Consequently, PCA is a recommended tool for analysis of such  high-dimensional data.

There is an underlying probabilistic model behind PCA, called probabilistic principal component analysis (PPCA), defined as follows. The observation  vectors  $\{\bx_i\}_{1\le i\le n}$ are
$p$-dimensional and  satisfy the equation
\begin{equation}\label{model}
  \bx_i=\bLa \bdf_i + \bde_i +\bmu ~, \quad i=1,\ldots,n.
\end{equation}
Here,
$\bdf_i$ is a $m$-dimensional {\em principal components} with $m\ll p$,
$\bLa$  is a $p\times m$ matrix of {\em loadings},   and
$\bmu$ represents the general mean
and $(\bde_i)$ are a sequence of independent Gaussian errors with covariance matrix $\bPsi=\sigma^2 \bI_p$. The parameter $\sigma^2$ is the noise variance we are interested in.
The components $(\bdf_i)$ are also Gaussian and unobserved.

To ensure the identification of the model, constraints have to be introduced on the parameters. There are
several possibilities for the choice of such constraints, see
e.g. Table 1 in \citet{BaiLi12}.
A traditional choice is the following \citep[Chapter 14]{Anderson3}:
\begin{itemize}
\item   $ \E \bdf_i= {\bf 0}$ and
  $\E\bdf_i\bdf_i'=\bI$;
\item
  The matrix $\bGa:=\bLa'\bLa$  is diagonal
  with distinct diagonal elements.
\end{itemize}
Therefore, the population covariance matrix (PCM) of $\{\bx_i\}$ is
\begin{equation}\label{eq:Sigma}
  \bSi=\bLa\bLa'+\si^2 \bI.
\end{equation}
Finding a reliable estimator of $\si^2$ is a nontrivial issue for high-dimensional data which we now pursue.

The PPCA model (1) can be viewed as a special instance of the approximate factor model \citep{Chamberlain} when the noise covariance $\bPsi$ is a general diagonal matrix. (The model is also called a strict factor model in statistical literature, see \citet{Anderson3}). For related recent papers on inference of large approximate (or dynamic) factor models, we refer to \citet{BaiJ03}, \citet{forni2000reference} and \citet{doz2012quasi}.

Let $\bar{\x}$ be the sample mean and define the sample covariance matrix
\begin{eqnarray}
\mathbf{S}_n= \frac{1}{n-1} \sum_{i=1}^n (\mathbf{x}_i - \bar{\x})
(\mathbf{x}_i -\bar{\x})'\text{.}
\end{eqnarray}
Let $\lambda_{n,1} \ge \lambda_{n,2} \ge \dots \ge
\lambda_{n,p}$ be the eigenvalues of $\mathbf{S}_n$. The maximum likelihood estimator of the noise variance is
\begin{eqnarray}
  \widehat{\sigma}^2 &=& \frac{1}{p-m} \sum_{i=m+1}^p \lambda_{n,i}.
\end{eqnarray}

In the classic setting where the dimension $p$ is relatively small compared to the sample size $n$ (low-dimensional setting), the consistency of $\widehat{\sigma}^2$ is established in \citet{Anderson1}. Moreover, it is asymptotically normal with the standard $\sqrt{n}$-convergence rate, as $n \rightarrow \infty$,
\begin{equation}
  \sqrt n ( \widehat{\sigma}^2 -{\sigma}^2 )
  \overset{\mathcal{D}}{\longrightarrow} \mathcal{N}(0, s^2), \quad
  s^2= \frac{2\sigma^4}{p-m}\text{.}\label{sigmaMLE}
\end{equation}
Actually,  \citet{Anderson2} provides a general CLT for m.l.e. in an approximate factor model that encompasses the present PPCA model. For the reader's convenience, we provide in Subsection 6.5 a detailed deviation of (5) from this general CLT.

The situation is, however radically different when $p$ is large compared to the sample size $n$. Recent advances in high-dimensional statistics indicate that in such high-dimensional situation, the above asymptotic result is no more valid and indeed, it has been widely observed in the literature that $\widehat{\sigma}^2$ seriously underestimates the true noise variance $\sigma^2$. Basically, when $p$ becomes large, the sample principal eigenvalues and principal components are no longer consistent estimates of their population counterparts \citep{baik2006eigenvalues, johnstone2009consistency, kritchman2008determining}. Many estimation methods developed in low-dimensional setting have been shown to perform poorly even for moderately large $p$ and $n$ \citep{cragg1997inferring}.

As all meaningful inference procedures in the model will unavoidably use some estimator of the noise variance $\sigma^2$, such a severe bias needs to be corrected for high-dimensional data. There are several estimators proposed to deal with the high-dimensional situation. \cite{Kritchman1} proposes an estimator by solving a system of implicit equations; \citet{Solo} introduces an estimator using the  median of the sample eigenvalues $\{\lambda_{n,k}\}$; and \citet{johnstone2009consistency} uses the median of the sample variances. However, these estimators are assessed by Monte-Carlo experiments only and their theoretical properties have not been investigated.

The main aim of this paper is to provide a new estimator  of the noise variance for which a rigorous asymptotic theory can be established in the high-dimensional setting. First, by using recent advances in random matrix theory, we found a CLT for the m.l.e. $\widehat{\sigma}^2$ in the high-dimensional setting (Theorem 1 in Section 2). Next, using this identification and random matrix theory, we propose a new estimator $\widehat{\sigma}_\ast^2$ for the noise variance by correcting this bias. The asymptotic normality of the new estimator is thus established with explicit asymptotic mean and variance.

Although the asymptotic Gaussian distribution of the bias-corrected estimator $\widehat{\sigma}_\ast^2$ is established under the high-dimensional setting $p\to \infty$, $n\to \infty$ and $p/n \to c >0$, if we set $c=0$, i.e. the dimension $p$ is infinitely smaller than $n$, this Gaussian limit coincides with the classical low-dimensional limit given in (5). In this sense, the new asymptotic theory extends in a continuous manner the classical low-dimensional result to the high-dimensional situation. Finite sample properties of the bias-corrected estimator $\widehat{\sigma}_\ast^2$ have been checked via Monte-Carlo experiments in comparison with the above-mentioned three existing estimators. In terms of mean squared errors and in all the tested scenarios, $\widehat{\sigma}_\ast^2$ outperforms very significantly two of them, and is slightly preferable than the third one, see Table 3.

In order to demonstrate further potential benefits achievable by the implementation of the bias-corrected estimator $\widehat{\sigma}_\ast^2$, we develop three applications. Two of these applications concern an important inference problem in PPCA, namely, the determination of the number of principal components (PCs). This problem involves the noise variance estimation and many methods for choosing PCs have been proposed in the literature. We consider two benchmark procedures from different fields. The first one is proposed in \citet{Solo} that employs PPCA to signal processing and proposes the Stein's unbiased risk estimator (SURE) to find the number of PCs. The second one is from \citet{Baij} which develops six criteria with penalty on both $p$ and $n$ to identify the number of factors in the approximate factor model. The approximate factor model allows the errors $(\bde_i)$ be correlated. PPCA can be considered as a simplified instance of this model and indeed, \citet{Baij} applies also their criteria to PPCA. Notice that the determination criteria in \citet{Baij} are popular and widely cited in the literature. Furthermore, both the procedures in \citet{Solo} and \citet{Baij} implement a specific estimator of the noise variance $\sigma^2$. By substituting the bias-corrected estimator $\widehat{\sigma}_\ast^2$ for their estimators, we demonstrate by extensive simulation that these procedures are uniformly and significantly improved.

The third application of the bias-corrected estimator $\widehat{\sigma}_\ast^2$ concerns the goodness-of-fit test for the PPCA model. The LR test statistic as well as their classical (low-dimensional) chi-squared asymptotics are well-known since the work of \citet{Amemiya}. These results are again challenged by high-dimensional data and the classical chi-squared limit is no more valid. Following an approach devised in \citet{Bai2}, we propose a correction to this goodness-of-fit test statistic involving our new estimator $\widehat{\sigma}_\ast^2$ of noise variance to cope with the high-dimensional effects and establish its asymptotic normality.

The remaining sections are organized as follows. In Section 2, we present the main results of the paper. In Section 3, our new estimator $\widehat{\sigma}_\ast^2$ is substituted for an estimator from the authors of the SURE criterion and it is shown that this substitution improves greatly the SURE criterion. Similarly, in Section 4, the estimator $\widehat{\sigma}_\ast^2$ is applied to the criteria proposed by \citet{Baij} to determine the number of PCs and again an improvement of these criteria is obtained. In Section 5, we develop the corrected likelihood ratio test for the goodness-of-fit of a PPCA model in the high-dimensional framework using the new estimator $\widehat{\sigma}_\ast^2$. Technical proofs are gathered together in Section 6. Section 7 concludes.

\section{Main results}
The PPCA model (1) is a spiked population model \citep{Johnstone} since the eigenvalues of PCM $\bSi$ are
\begin{eqnarray}
  \mbox{spec}(\bSi) &=&  \nonumber
  ( \underbrace{\alpha_1,\dots,\alpha_1}_{n_1},\dots,\underbrace{\alpha_K,\dots,\alpha_K}_{n_K},\underbrace{0,\dots,0}_{p-m}  ) + \sigma^2 (\underbrace{1,\dots,1}_p)\\
&=&\sigma^2( \underbrace{\alpha_1^*,\dots,\alpha_1^*}_{n_1},\dots,\underbrace{\alpha_K^*,\dots,\alpha_K^*}_{n_K},\underbrace{1,\cdots,1}_{p-m}  )\text{,} \label{6spike_model}
\end{eqnarray}
where $(\alpha_i)$ are non-null eigenvalues of $\bLa\bLa'$ with
multiplicity numbers $(n_i)$ satisfying $n_1+\cdots+n_K=m$ and the notation $\alpha_i^* = {\alpha_i}{/\sigma^2} + 1$ is used. To develop a meaningful asymptotic theory in the high-dimensional context, we assume that $p$ and $n$ are related so that when $n
\rightarrow \infty$, $c_n=p/(n-1) \rightarrow c > 0$, that is, $p$
can be large compared to the sample size $n$ and for the asymptotic theory,  $p$ and $n$ tend to infinity  proportionally. Define the function
\[\phi(\alpha) = \alpha + \frac{c\alpha}{\alpha-1}~\text{,}\quad
\alpha\neq 1~,\]
and set $s_0=0$ and $s_i=n_1+\cdots+n_i$
for $1\le i\le K$. The set
$J_i=\{s_{i-1}+1,\dots,s_i\}$ is then the indexes  among
$\{1,\ldots,p\} $ associated to $\alpha_i$ counting the
multiplicities.
Following   \citet{Baik},
assumed that
$\alpha_1^* \ge \dots \ge \alpha_{m}^* > 1+\sqrt{c}$, i.e all the
eigenvalues $\alpha_i$ are greater than $\sigma^2\sqrt{c}$.
It is then known that, for the spiked sample eigenvalues $\lambda_{n,k}$ of $\mathbf{S}_n$,
$1\le k\le m$,  almost surely if $k\in J_i$,
\begin{eqnarray}
  \lambda_{n,k} &\longrightarrow& \sigma^2\phi(\alpha_i^*)
  =  \alpha_i + \sigma^2 + \sigma^2c \left ( 1 + \frac{\sigma^2}{\alpha_i}\right ) \label{llambda} \text{.}
\end{eqnarray}
Moreover, the remaining sample eigenvalues  $\{\lambda_{n,k}\}_{m<k\le
p}$, called {\em noise eigenvalues}, will converge to a continuous
distribution  with support interval $[a(c), b(c)]$
where $a(c)=\sigma^2(1-\sqrt{c})^2$
and $b(c)=\sigma^2(1+\sqrt{c})^2$.
In particular,
for all $1\le i \le L$ with a prefixed range $L$ and almost surely,
$\lambda_{n,m+i} \rightarrow b(c)~.$ It is worth noticing
 that in \eqref{llambda}, if we let $c\to 0$, we recover the
low-dimensional limits  $\lambda_{n,k}\to \alpha_i+\sigma^2$ (population spike
eigenvalues)
and  $\lambda_{n,k}\to \sigma^2$ (population noise eigenvalues) discussed
earlier.

In a further step, CLT for the spiked eigenvalues is established in
\citet{Bai1} (see also \citet{Paul}): the $n_i$-dimensional vector
\begin{eqnarray}
  & \{\sqrt{n}(\lambda_{n,k} - \sigma^2\phi(\alpha_k^*))\mbox{,} k \in J_i \}& \label{dist}
\end{eqnarray}
converges to a well-determined
$n_i$-dimensional limiting distribution. Moreover,
this limiting distribution
is  Gaussian if and only if the
corresponding population spike eigenvalue $\alpha_i$ is simple,
i.e. $n_i=1$.

As explained in the Introduction, when the dimension $p$ is large compared to the sample size $n$, the m.l.e. $\widehat{\sigma}^2$ in (4) has a negative bias. In order to identify  this bias, we first
establish a central limit theorem for $\widehat{\sigma}^2$ under the high-dimensional scheme.
\begin{theorem2} \label{T1}
 Consider the PPCA model (1) with population
 covariance matrix $\bSi=\bLa\bLa'+\sigma^2 \I_p$
 where both the principal components and the noise are Gaussian.
 Assume that $p\to\infty$,
 $n\to\infty$ and $c_n=p/(n-1)\to c>0$. Then, we have
  \[\frac{(p-m)}{\sigma^2\sqrt{2c}} (\widehat{\sigma}^2 - \sigma^2) + b(\sigma^2)\overset{\mathcal{D}}{\longrightarrow} \mathcal{N}(0,1)\text{,}\]
  where
  $b(\sigma^2)=\sqrt{\frac{c}{2}} \left ( m + \sigma^2 \sum_{i=1}^m \frac{1}{\alpha_i}\right )$.
\end{theorem2}

Therefore, for high-dimensional data, the m.l.e. $\widehat{\sigma}^2$ has an asymptotic  bias
$ - b(\sigma^2) $ (after normalization). This bias is a complex
function of the noise variance  and the $m$ non-null eigenvalues of the loading
matrix $\bLa\bLa'$.  Notice that  the above CLT is
still valid if    $\tilde c_n=(p-m)/n$  is substituted for $c$. Now if indeed $p\ll n$, i.e. the dimension $p$ is infinitely smaller than the sample size $n$, so that  $\tilde c_n \simeq  0$ and $b( \sigma^2)\simeq
0$, and hence
\[ \frac{(p-m)}{\sigma^2\sqrt{2c}} (\widehat{\sigma}^2 - \sigma^2) +
b(\sigma^2) \simeq  \frac{\sqrt{p-m}}{\sigma^2\sqrt{2}}    \sqrt{n}(\widehat{\sigma}^2 - \sigma^2)\overset{\mathcal{D}}{\longrightarrow} \mathcal{N}(0,1)~.
\]
This is nothing but the CLT \eqref{sigmaMLE} for $\widehat{\sigma}^2$ known under  the classical low-dimensional scheme. In a sense, Theorem~\ref{T1} constitutes a natural and continuous extension of the classical CLT to the high-dimensional
context.

The previous theory recommends to correct the negative bias of
$\widehat{\sigma}^2$. However, the bias $b(\sigma^2)$ depends on the number $m$ and the values $\alpha_i$ of the spikes. These parameters are likely unknown in real-data applications and they need to be estimated. In the literature, consistent estimators of $m$ have been proposed, e.g. in \citet{Baij}, \citet{Kritchman1}, \citet{Onatski09} and \cite{Passemier1}. For the values of the spikes $\alpha_i$, it is
easy to see that it can be estimated by inverting the function $\phi$ in
(\ref{llambda}) at the corresponding eigenvalues
$\lambda_{n, j}$. Moreover, by applying the delta-method to (\ref{dist}),
we can obtain the asymptotic distribution of this estimator, see
\citet{BaiDing13}.

As the bias depends on $\sigma^2$ which we want to estimate, a natural correction is to use the plug-in estimator
\begin{eqnarray}
\widehat{\sigma}_*^2=\widehat{\sigma}^2+\frac{b(\widehat{\sigma}^2)}{p-m}\widehat{\sigma}^2 \sqrt{2c_n}\text{.}
\end{eqnarray}
This estimator will be hereafter referred as the {\em bias-corrected estimator}. Notice that here the number of factors $m$ can be
replaced by any consistent estimate as discussed above   without affecting the limiting distribution of the estimator. The following CLT is an easy consequence of Theorem 1.

\begin{theorem2}
  We assume the same conditions as in Theorem 1.
  Then, we have
  \[  \frac{p-m}{\sigma^2\sqrt{2c_n}}    \left(   \widehat{\sigma}_*^2 - \sigma^2     \right) \overset{\mathcal{D}}{\longrightarrow} \mathcal{N}(0,1)~.
\]
\end{theorem2}
Compared to the m.l.e. $\widehat{\sigma}^2$ in Theorem~\ref{T1}, the bias-corrected estimator $\widehat{\sigma}_\ast^2$ has no more a bias after normalization by $\frac{p-m}{\sigma^2\sqrt{2c_n}}$, and it should be a much better estimator than $\widehat{\sigma}^2$.

\subsection{Monte-Carlo experiments}
We first check by simulation the effect of bias-correction obtained in $\widehat{\sigma}_\ast^2$ and its asymptotic normality. Independent Gaussian samples of size $n$ are considered in three different settings:
\begin{list}{$\bullet$}{\leftmargin=2em}
\item Model 1: $\mbox{spec}(\bSi)= (25,16,9,0,\dots,0)+\sigma^2(1,\dots,1)$, $\sigma^2=4$, $c=1$;
\item Model 2: $\mbox{spec}(\bSi)= (4,3,0,\dots,0)+\sigma^2(1,\dots,1)$, $\sigma^2=2$, $c=0.2$;
\item Model 3: $\mbox{spec}(\bSi)= (12,10,8,8,0,\dots,0)+\sigma^2(1,\dots,1)$, $\sigma^2=3$, $c=1.5$.
\end{list}

In Table \ref{t1c6}, we compare the empirical bias of
$\widehat{\sigma}^2$ (i.e. the empirical mean of
$\widehat{\sigma}^2-\sigma^2=\frac{1}{p-m} \sum_{i=m+1}^p
\lambda_{n,i}-\sigma^2$) over 1000 replications with the theoretical  one $-\sigma^2\sqrt{2c}  b(\sigma^2)/(p-m)$ in different settings.
In all the three models, the empirical and theoretical bias are  close each other. As expected, their  difference vanishes when $p$ and $n$ increase. The table also shows that this bias is quite significant even for large dimension and sample size such as $(p, n)=(1500, 1000)$.

{
\begin{table}[!ht]
\begin{center}
\caption{{  Comparison between the empirical and the theoretical bias in various settings.}} \label{t1c6}
\small
{ \begin{tabular}{|cll|ccc|}
\hline
\multicolumn{3}{|c|}{Settings}  & Empirical bias & Theoretical bias & $|$Difference$|$ \\
\hline
\multirow{3}{*}{Model 1} & $p=100$ & $n=100$ & -0.1556 & -0.1589 & 0.0023\\
 & $p=400$ & $n=400$  & -0.0379 & -0.0388 & 0.0009\\
 & $p=800$ & $n=800$  & -0.0189 & -0.0193 & 0.0004\\
\hline
\multirow{3}{*}{Model 2} & $p=20$ & $n=100$ & -0.0654 & -0.0704 & 0.0050\\
 & $p=80$ & $n=400$ & -0.0150 & -0.0162 & 0.0012\\
 & $p=200$ & $n=1000$  & -0.0064 & -0.0063 & 0.0001\\
\hline
\multirow{3}{*}{Model 3} & $p=150$ & $n=100$ & -0.0801 & -0.0795 & 0.0006\\
 & $p=600$ & $n=400$  & -0.0400 & -0.0397 & 0.0003\\
 & $p=1500$ & $n=1000$  & -0.0157 & -0.0159 & 0.0002\\
\hline
\end{tabular}}
\end{center}
\end{table}}

\begin{figure}
\centering
   \subfigure[]{
\includegraphics[width=0.42\columnwidth]{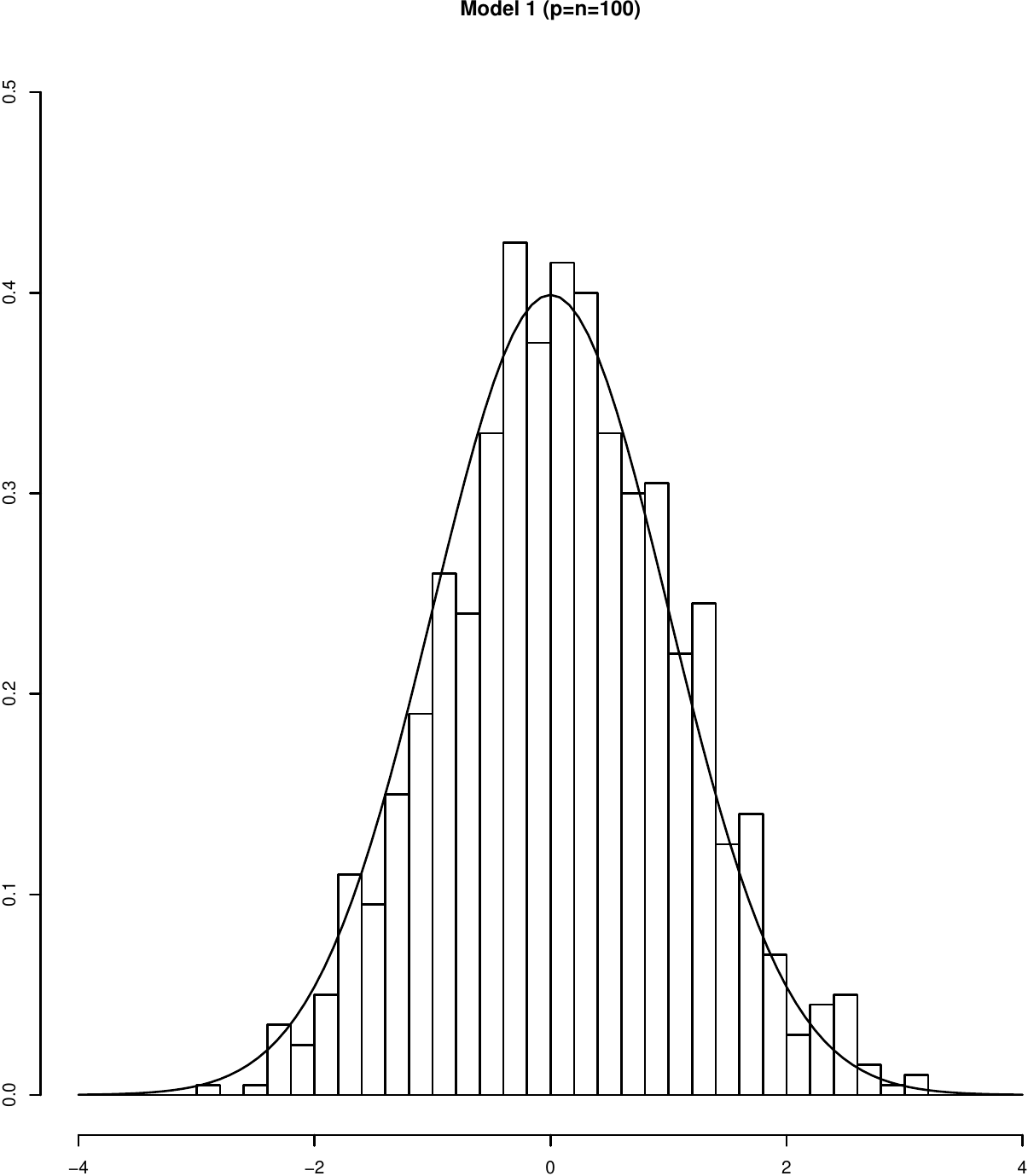}}
\subfigure[]{
\includegraphics[width=0.42\columnwidth]{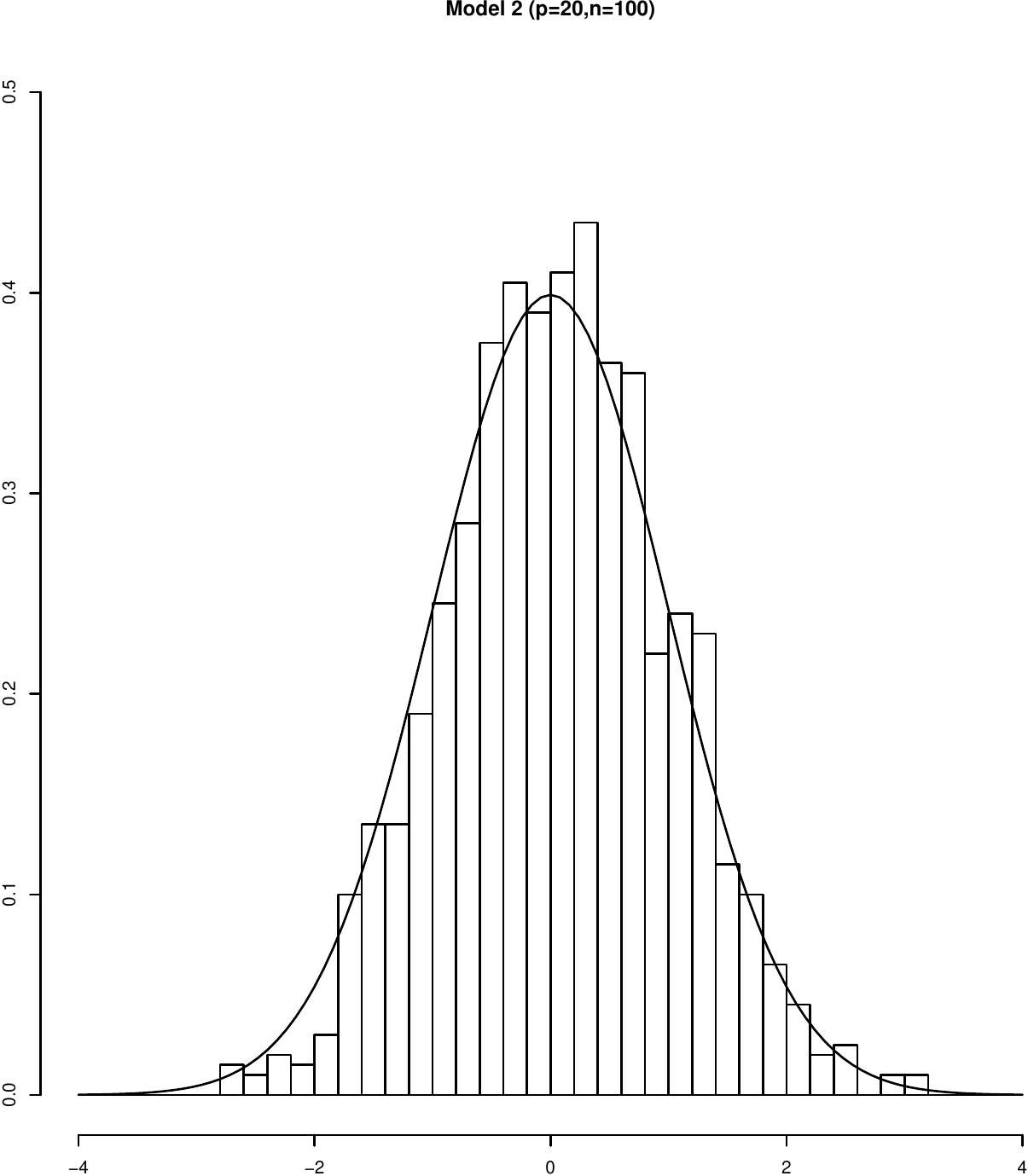}}\\
\subfigure[]{
\includegraphics[width=0.42\columnwidth]{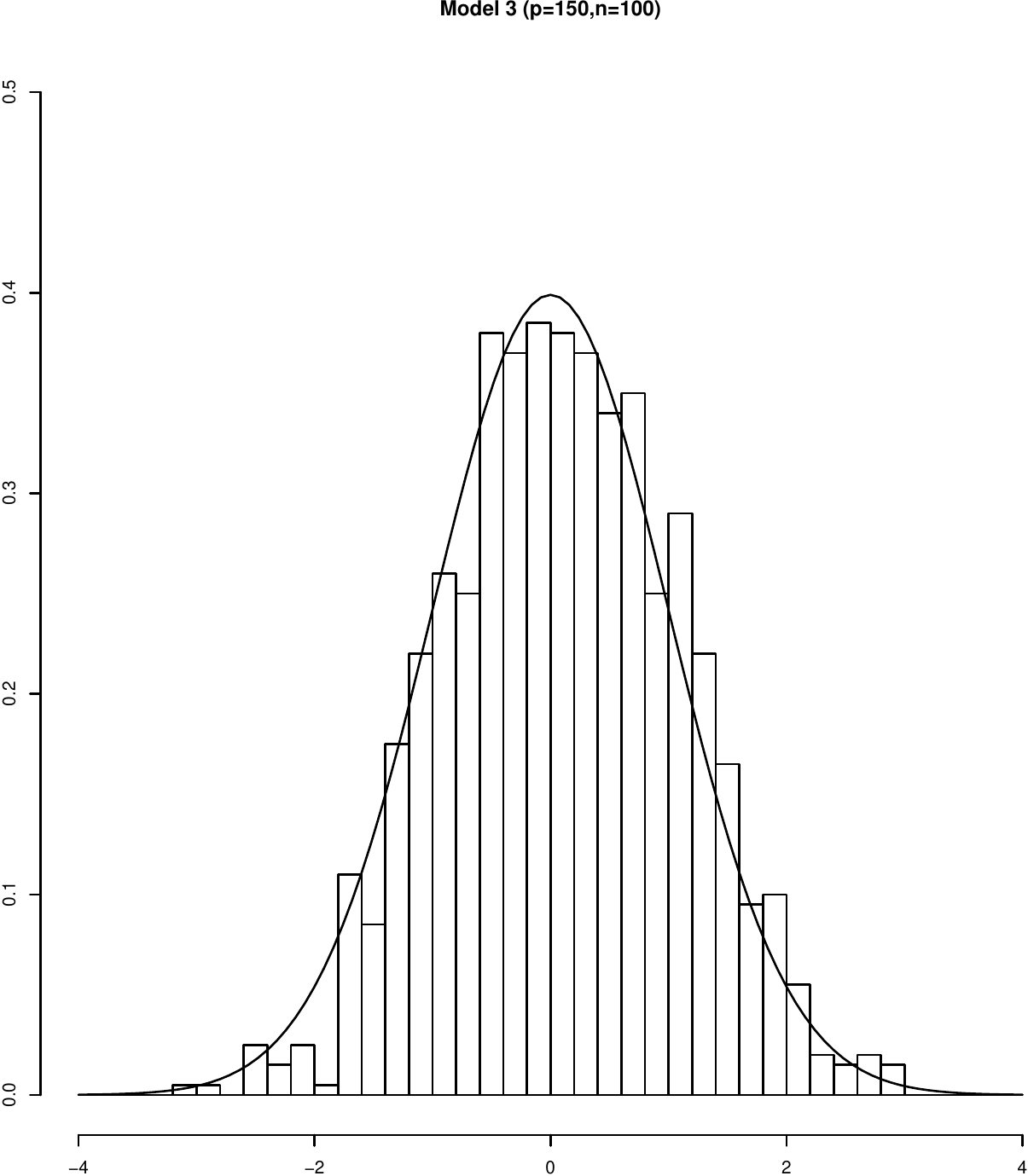}}
   \caption{{\footnotesize Histogram of $\frac{(p-m)}{\sigma^2\sqrt{2c}} (\widehat{\sigma}^2 - \sigma^2) + b(\sigma^2)$ compared with the density of a standard Gaussian distribution.}} \label{f1c6}
\end{figure}

Figure \ref{f1c6} presents the histograms from  1000 replications of
\[\frac{(p-m)}{\sigma^2\sqrt{2c_n}} (\widehat{\sigma}^2 - \sigma^2) + b(\sigma^2)\]
of the three models above, with sample size $n=100$ and dimensions $p=c
\times n$, compared to the density of the standard Gaussian distribution. The sampling distribution is almost normal.

To assess the quality of the bias-corrected estimator
$\widehat{\sigma}_*^2$, we conduct some simulation experiments using the same settings. Table \ref{t2c6} compares
the empirical means, mean absolute deviations (MAD) and the mean squared errors (MSE) of $\widehat{\sigma}^2$ and $\widehat{\sigma}_*^2$  over 1000 replications. The improvement of $\widehat{\sigma}_*^2$ over $\widehat{\sigma}^2$ is very significant in term of MAD in all the tested cases (the difference in term of MSE is less remarkable because the squares become quite small numbers).

Next, we compare our bias-corrected estimator to three existing estimators in the literature. For the reader's convenience, we recall their definitions:
\begin{enumerate}
\item The estimator  $\widehat{\sigma}_{\text{KN}}^2$ of \cite{Kritchman1}: it is defined as the solution of the following non-linear system of $m+1$ equations involving the $m+1$ unknowns $\widehat{\rho}_1,\dots,\widehat{\rho}_m$ and $\widehat{\sigma}_{\text{KN}}^2$:
  \begin{align*}
    \widehat{\sigma}_{\text{KN}}^2  - \frac{1}{p-m} \left [ \sum_{j=m+1}^p \lambda_{n,j} + \sum_{j=1}^m (\lambda_{n,j}- \widehat{\rho}_j) \right ] &=0\text{,} \quad \text{and}\\
    \widehat{\rho}_j^2-\widehat{\rho}_j\left ( \lambda_{n,j} + \widehat{\sigma}_{\text{KN}}^2 - \widehat{\sigma}_{\text{KN}}^2 \frac{p-m}{n}\right ) +\lambda_{n,j}\widehat{\sigma}_{\text{KN}}^2 &= 0\text{,} \quad j=1, \ldots, m.
  \end{align*}
  We use the computing code  available on the authors' web-page to carry
  out the simulation.
\item The estimator $\widehat{\sigma}_{\text{US}}^2$ of \cite{Solo}: it is defined as the ratio
  \[\widehat{\sigma}_{\text{US}}^2= \frac{\text{median}(\lambda_{n,m+1},\dots, \lambda_{n,p})}{m_{p/n,1}}\text{,}\]
  where $m_{\alpha,1}$ is the median of the Mar{\v c}enko-Pastur
  distribution $F_{\alpha,1}$ (more details on this distribution is given in Subsection 6.1).
\item The estimator $\widehat{\sigma}_{\text{median}}^2$ of \cite{johnstone2009consistency}: it is defined as the median of the $p$ sample variances
    \[\widehat{\sigma}_{\text{median}}^2=\text{median}\left(\frac{1}{n}\sum_{i=1}^n x_{ij}^2, \quad 1\leq j \leq p\right) \text{.} \]
    Clearly for this estimator, the data $\{x_{ij}\}$ are assumed centered.
\end{enumerate}

{
\begin{table}[!ht]
\caption{{ Empirical mean, MAD and MSE of $\widehat{\sigma}^2$ and $\widehat{\sigma}_*^2$ in various settings.}} \label{t2c6}
\centering
\small
{\begin{tabular}{|cccc|lcr|lcr|}
\hline
\multicolumn{4}{|c|}{Settings}  & \multicolumn{1}{c}{\multirow{2}{*}{$\widehat{\sigma}^2$}} & \multirow{2}{*}{MAD} & \multirow{2}{*}{MSE} & \multicolumn{1}{c}{\multirow{2}{*}{$\widehat{\sigma}_*^2$}} & \multirow{2}{*}{MAD} &\multirow{2}{*}{MSE}\\
Mod. & p & n & $\sigma^2$ & & & & & &\\
\hline
\multirow{3}{*}{1} & $100$ & $100$ & \multirow{3}{*}{4} & 3.8464 &0.1536 & 0.0032 & 3.9979 & 0.0021 & 0.0035\\
 & $400$ & $400$  &  & 3.9616 &0.0384& 0.0002 & 4.0000 & $< 10^{-5}$ & 0.0002\\
 & $800$ & $800$  &  & 3.9809 &0.0191& 0.0001 & 4.0002 &0.0002 & 0.0001\\
\hline
\multirow{3}{*}{2} & $20$ & $100$ & \multirow{3}{*}{2} &  1.9340 &0.0660 & 0.0043 & 2.0012 &0.0012& 0.0047\\
 & $80$ & $400$  &  & 1.9841 &0.0159& 0.0003 & 2.0001 &0.0001& 0.0003\\
 & $200$ & $1000$  &  & 1.9939 &0.0061&$< 10^{-5}$ & 2.0002 &0.0002& $< 10^{-5}$\\
\hline
\multirow{3}{*}{3} & $150$ & $100$ & \multirow{3}{*}{3} & 2.8400&0.1600 & 0.0011 & 2.9926 &0.0074& 0.0013\\
 & $600$ & $400$  &  & 2.9605 & 0.0395&0.0001 & 2.9999 &0.0001& 0.0001\\
 & $1500$ & $1000$  &  & 2.9839& 0.0161 & $< 10^{-5}$ & 2.9998 &0.0002 & $< 10^{-5}$\\
\hline
\end{tabular}}
\end{table}}

\begin{table}
\caption{Comparison between three existing estimators and the proposed one in terms of ratios of MSE: $\frac{\text{MSE}\big(\widehat{\sigma}_{\text{KN}}^2 \big)}{\text{MSE}\big(\widehat{\sigma}_\ast^2\big)}$, $\frac{\text{MSE}\big(\widehat{\sigma}_{\text{US}}^2 \big)}{\text{MSE}\big(\widehat{\sigma}_\ast^2\big)}$ and $\frac{\text{MSE}\big(\widehat{\sigma}_{\text{median}}^2 \big)}{\text{MSE}\big(\widehat{\sigma}_\ast^2\big)}$.}
\centering
\small
\begin{tabular}{|cccc|ccc|}
\hline
\multicolumn{4}{|c|}{Settings} & \multicolumn{1}{c}{\multirow{2}{*}{$\widehat{\sigma}_{KN}^2$}} & \multirow{2}{*}{$\widehat{\sigma}_{US}^2$} & \multirow{2}{*}{$\widehat{\sigma}_{\text{median}}^2$}\\
Mod. & p & n & $\sigma^2$ & & &  \\
\hline
\multirow{3}{*}{1} & $100$ & $100$ & \multirow{3}{*}{4}  & 1.01 & 4.40 & 1.47\\
& $400$ & $400$  &  & 1.00 & 6.50 & 1.59\\
 & $800$ & $800$  &  & 1.00 & 4.00 & 1.62\\
  \hline
\multirow{3}{*}{2} & $20$ & $100$ & \multirow{3}{*}{2} & 1.04 & 1.85 & 1.67\\
 & $80$ & $400$  &  &  1.00 & 2.67 & 1.52\\
 & $200$ & $1000$  &  &  1.00 & 10.10 & 1.53\\
  \hline
\multirow{3}{*}{3} & $150$ & $100$ & \multirow{3}{*}{3} & 1.07 & 7.08 & 1.26\\
 & $600$ & $400$  &  & 1.00 & 7.00 & 1.52\\
 & $1500$ & $1000$  &  & 0.96 & 10.10 & 1.60\\
\hline
\end{tabular}
\end{table}

Table 3 gives the ratios of the empirical MSEs of these estimators over the empirical MSE of the bias-corrected estimator $\widehat{\sigma}^2_*$.
The performances of $\widehat{\sigma}^2_*$ and
$\widehat{\sigma}^2_{\text{KN}}$ are similar and $\widehat{\sigma}^2_*$ being slightly better.
The estimator $\widehat{\sigma}^2_{\text{US}}$ shows slightly better
than the m.l.e. $\widehat{\sigma}^2$. The estimator $\widehat{\sigma}_{\text{median}}^2$ is better than $\widehat{\sigma}^2_{\text{US}}$ and the m.l.e. $\widehat{\sigma}^2$.
But $\widehat{\sigma}_{\text{median}}^2$ and $\widehat{\sigma}^2_{\text{US}}$ performs poorly
compared to $\widehat{\sigma}^2_*$ and
$\widehat{\sigma}^2_{\text{KN}}$.
The reader is, however reminded that the theoretic properties of  $\widehat{\sigma}^2_{\text{KN}}$,
$\widehat{\sigma}^2_{\text{US}}$ and $\widehat{\sigma}_{\text{median}}^2$ are unknown and so far there have been checked via simulations only.

\section{Application to the SURE criterion}
\citet{Solo} proposes to use Stein's unbiased risk estimator, SURE, to choose the number of PCs. This criterion uses the noise variance estimator $\widehat{\sigma}_{US}^2$ defined in Section 2. It aims at minimizing the Euclidean distance between the underlying estimator of the population mean $\bmu$ and its true value. The proposed SURE criterion for $m$ number of PCs (to be minimized) is
\begin{eqnarray}
R_m&=&(p-m)\widehat{\sigma}_{US}^2+\widehat{\sigma}_{US}^4\sum_{j=1}^m \frac{1}{\lambda_j} +2\widehat{\sigma}_{US}^2 (1-1/n)m \nonumber\\
&&-2\widehat{\sigma}_{US}^4(1-1/n)\sum_{j=1}^m \frac{1}{\lambda_j} + \frac{4(1-1/n)\widehat{\sigma}_{US}^4}{n}\sum_{j=1}^m \frac{1}{\lambda_j} +C_m,
\end{eqnarray}
where
\begin{eqnarray*}
C_m&=&\frac{4(1-1/n)\widehat{\sigma}_{US}^2}{n}\sum_{j=1}^m \sum_{i=m+1}^p \frac{\lambda_j-\widehat{\sigma}_{US}^2}{\lambda_j-\lambda_i} + \frac{2(1-1/n)\widehat{\sigma}_{US}^2}{n}m(m-1)\\
&& \quad - \frac{2(1-1/n)\widehat{\sigma}_{US}^2}{n}(p-1)\sum_{j=1}^m \left(1-\frac{\widehat{\sigma}_{US}^2}{\lambda_j}\right).
\end{eqnarray*}
Recall that $\widehat{\sigma}_{US}^2$ is also related to $m$. From Section 2, we have known that $\widehat{\sigma}_{US}^2$ is not as good as our bias-corrected estimator. To examine further this difference, we replace  $\widehat{\sigma}_{US}^2$ with $\widehat{\sigma}_\ast^2$ in (10), referred then as SURE$^\ast$, to see whether the performance of SURE can be improved.

Then simulation experiments are conducted to compare SURE with SURE$^\ast$. The setup follows the paper \citet{Solo} and the data are simulated according to (1) with the parameters $p=64, p/n=[2/3, 1/2, 2/5], m=[5, 10, 15, 20]$ and $\sigma^2=1$.
The loading matrix is set to $\bLa=\mathbf{F}\mathbf{D}^{1/2}$, where $\mathbf{F}$ is constructed by generating a $p\times m$ matrix of Gaussian random variables and then orthogonalizing the resulting matrix, and
\begin{eqnarray*}
\mathbf{D} &=& \text{diag}(\lambda_1, \lambda_2, \ldots, \lambda_{m-1}, \lambda_m)\\
&=&\text{diag}\big((m+1)^2, m^2, \ldots, 3^2, \lambda_m\big), \quad \lambda_m=1.5.
\end{eqnarray*}

All simulations were repeated $1500$ times. Table 4 shows the percentage of correct selection of PCs for SURE and SURE$^\ast$. It can be seen that SURE$^\ast$ largely outperforms SURE in all of the tested cases, most of times by a wide margin. All the percentages of correct selection of SURE$^\ast$ are larger than 90\% and in 4 out of 12 cases, the detection rate is 100\%. Therefore, by implementing our bias-corrected estimator of the noise variance instead of the one provided by its authors, the criterion SURE has a much better performance.
\begin{table}
\caption{Comparison between SURE and SURE$^\ast$ in terms of percentage of correct selection of PCs.}
\centering
\small
\begin{tabular}{|c|c|ccc|}
\hline
m & \ &  $p/n=2/3$ & $p/n=1/2$ & $p/n=2/5$\\
\hline
\multirow{2}{*}{5}  & SURE$^\ast$ & 1.000 & 1.000 & 1.000\\
\ & SURE & 0.408 & 0.621 & 0.807\\
\hline
\multirow{2}{*}{10} & SURE$^\ast$ & 0.990 & 1.000 & 0.998\\
\ & SURE & 0.512 & 0.739 & 0.858\\
\hline
\multirow{2}{*}{15} & SURE$^\ast$ &  0.904 & 0.978 & 0.989\\
\ & SURE & 0.598 & 0.783 & 0.911\\
\hline
\multirow{2}{*}{20} & SURE$^\ast$ & 0.908 &0.966 & 0.990\\
\ & SURE & 0.617 & 0.810 & 0.899\\
\hline
\end{tabular}
\end{table}

\section{Application to determination criteria of \citet{Baij}}
In econometrics, the assumption of additive white noise is reasonable for asset returns with low frequencies (e.g., monthly returns of stocks) \citep{ma2013sparse}. \citet{Baij} proposes six criteria to estimate the number of PCs (or factors) under the framework of large cross-sections ($N$) and large time dimensions ($T$). These criteria penalize both dimensions $N$ and $T$ and outperform the usual AIC and BIC, which are functions of $N$ or $T$ alone, under the assumption that both $N$ and $T$ grow to infinity. Notice that the dimension-sample-size pair is denoted here as $(N, T)$ instead of $(p, n)$. These six criteria are as follows:
\begin{eqnarray}
&&PC_{p1}(m)=V(m,\hat{F}^m) +m \widehat{\sigma}_{\text{BN}}^2\left(\frac{N+T}{NT}\right) \ln\left(\frac{NT}{N+T}\right);\nonumber\\
&&PC_{p2}(m)=V(m,\hat{F}^m) +m\widehat{\sigma}_{\text{BN}}^2\left(\frac{N+T}{NT}\right) \ln C_{NT}^2;\nonumber\\
&&PC_{p3}(m)=V(m,\hat{F}^m)+m\widehat{\sigma}_{\text{BN}}^2\left(\frac{\ln C_{NT}^2}{C_{NT}^2}\right);\nonumber\\
&&IC_{p1}(m)=\ln\left(V(m,\hat{F}^m)\right) +m \left(\frac{N+T}{NT}\right) \ln\left(\frac{NT}{N+T}\right);\nonumber\\
&&IC_{p2}(m)=\ln \left(V(m,\hat{F^m})\right) +m\left(\frac{N+T}{NT}\right) \ln C_{NT}^2;\nonumber\\
&&IC_{p3}(m)=\ln \left(V(m,\hat{F}^m)\right)+m\left(\frac{\ln C_{NT}^2}{C_{NT}^2}\right),
\end{eqnarray}
where $\widehat{\sigma}_{\text{BN}}^2$ is a consistent estimate of $(NT)^{-1}\sum_{i=1}^N\sum_{j=1}^T E(e_{ij})^2$, $V(m,\hat{F}^m)=(NT)^{-1}\sum_{i=1}^N \hat{\bde}_i^\prime \hat{\bde}_i$, and $C_{NT}^2=\min\{N, T\}$. Notice that the difference between the $PC_p$ and the $IC_p$ criteria is that the first family uses the $V$ function while the second family uses its logarithm. In applications, $\widehat{\sigma}_{\text{BN}}^2$ is replaced by $V(m_0,\hat{F}^{m_0})$, where $m_0$ is a predetermined maximum value of $m$. We can see that the calculations of $\widehat{\sigma}_{\text{BN}}^2$ and $V(m,\hat{F}^m)$ have no explicit formula and are based on the estimation of the residuals $(\hat{\bde}_i)$. It is worth mentioning that $\widehat{\sigma}_{\text{BN}}^2$ is indeed the estimator of the noise variance if the underlying model is the PPCA model. Now we substitute the proposed bias-corrected estimator $\widehat{\sigma}_\ast^2$ for empirical $\widehat{\sigma}_{\text{BN}}^2$ and update in accordance the statistic $V(m,\hat{F}^m)$. The modified criteria using $\widehat{\sigma}_\ast^2$ are denoted as $PC^\ast_{p1}, PC^\ast_{p2}, PC^\ast_{p3}, IC^\ast_{p1}, IC^\ast_{p2}$ and $IC^\ast_{p3}$.

Simulation experiments are conducted to check whether the performance of the three $PC_p$ criteria is improved by the bias-corrected estimator. The results of three $IC_p$ criteria are not presented here since they are very similar in our context to the results of the $PC_p$ criteria. As in \citet{Baij}, the data are generated from the model:
\begin{eqnarray}
X_{it}=\sum_{j=1}^m \lambda_{ij}F_{tj} +\sqrt{\theta}e_{it},
\end{eqnarray}
where the PCs, the loadings and the errors $(e_{it})$ are $N(0,1)$ variates, the common component of $X_{it}$ has variance $m$ and the idiosyncratic component has variance $\theta$. Notice that the noise variance here is $\sigma^2=\theta$ and $\bLa=(\lambda_{ij})$. Typically, a PC corresponding to $\alpha_j$ is detectable when $\alpha_j\geq \left(1+\sqrt{\frac{N}{T}}\right)\theta$, see (7). We conduct extensive simulation by reproducing all the configuration of $N$ and $T$ used in \citet{Baij}. In particular, the last five rows of each table below correspond to small dimensions (either $N$ or $T$ is small).

Tables 5-8 report the empirical means of the estimator of the number of PCs over 1000 replications, for $m=1, 3$ and 5 respectively, with standard errors in parentheses. When a standard error is actually zero, no standard error is thus indicated. For all cases, the predetermined maximum number $m_0$ of PCs is set to 8. When the true number of PCs is 1 (Table 5), the criteria $PC_p^\ast$ can correctly detect the number almost surely and the corresponding standard errors are all zeros. In comparison, there are 11 cases where the criteria $PC_p$ lose efficiency in finding the true number of PCs with a non-zero standard error. In the small dimensions situation (last five rows),  all $PC_p^\ast$ and $PC_p$ fail when the value of $N$ is 10: they all report the maximum value $m_0$. But the criteria $PC_p^\ast$ outperform $PC_p$ in the last three cases in terms of mean and standard error. In Table 6, although the criteria $PC_p^\ast$ have some determination errors in a few cases, these errors are much smaller than the corresponding ones from the initial criteria $PC_p$. For the last three cases of small dimensions, the results of criteria $PC_p^\ast$ are also much better than that of $PC_p$. In Table 7, the common component and idiosyncratic component have the same variance 5 which is large and it can be verified that in this setting, some of the 5 PC eigenvalues $\alpha_j$ do not satisfy the detection condition $\alpha_j\geq \left(1+\sqrt{\frac{N}{T}}\right)\theta$. Consequently, the criteria $PC_p^\ast$ are no longer uniformly better than $PC_p$. In Table 8, when the variance of idiosyncratic component is smaller than that of the common component, the criteria $PC_p^\ast$ have again an overall better performance than the criteria $PC_p$. In conclusion, except some cases in Table 7, these criteria proposed in \cite{Baij} can have a better performance using the bias-corrected estimator proposed in this paper for PPCA models.

\begin{table}
\caption{Comparison between $PC^\ast_{p1}, PC^\ast_{p2}, PC^\ast_{p3}$ and $PC_{p1}, PC_{p2}, PC_{p3}$ in terms of  the mean estimation numbers of PCs for $m=1, \theta=1$.}
\centering
\tiny
\begin{tabular}{cc|ccc|ccc}
\hline
 N& T & $PC^\ast_{p1}$ & $PC^\ast_{p2}$ & $PC^\ast_{p3}$ & $PC_{p1}$ & $PC_{p2}$ & $PC_{p3}$\\
\hline
100 & 40 & 1.00 & 1.00 & 1.00 & 1.17(0.37) & 1.01(0.10) & 3.78(0.75)\\
100 & 60 & 1.00 & 1.00 & 1.00 & 1.00 & 1.00 & 3.63(0.76)\\
200 & 60 & 1.00 & 1.00 & 1.00 & 1.00 & 1.00 &1.00\\
500 & 60 &1.00 & 1.00 & 1.00 & 1.00 & 1.00 &1.00\\
1000 & 60 &1.00 & 1.00 & 1.00 & 1.00 & 1.00 &1.00\\
2000&60&1.00 & 1.00 & 1.00 & 1.00 & 1.00 &1.00\\
100 & 100 & 1.00 & 1.00 & 1.00 &1.00 & 1.00 & 5.36(0.80)\\
200 & 100 &1.00 & 1.00 & 1.00 & 1.00 & 1.00 &1.00\\
500 & 100 &1.00 & 1.00 & 1.00 & 1.00 & 1.00 &1.00\\
1000&100&1.00 & 1.00 & 1.00 & 1.00 & 1.00 &1.00\\
2000&100&1.00 & 1.00 & 1.00 & 1.00 & 1.00 &1.00\\
40 & 100 & 1.00 & 1.00 & 1.00 & 1.79(0.72) & 1.19(0.40) & 4.91(0.90)\\
60 & 100 & 1.00 & 1.00 & 1.00 & 1.01(0.08) & 1.00 & 4.30(0.85)\\
60 & 200 &1.00 & 1.00 & 1.00 & 1.00 & 1.00 &1.02(0.16)\\
60 &500 &1.00 & 1.00 & 1.00 & 1.00 & 1.00 &1.00\\
60 & 1000&1.00 & 1.00 & 1.00 & 1.00 & 1.00 &1.00\\
60&2000&1.00 & 1.00 & 1.00 & 1.00 & 1.00 &1.00\\
4000&60 &1.00 & 1.00 & 1.00 & 1.00 & 1.00 &1.00\\
4000&100&1.00 & 1.00 & 1.00 & 1.00 & 1.00 &1.00\\
8000&60&1.00 & 1.00 & 1.00 & 1.00 & 1.00 &1.00\\
8000&100&1.00 & 1.00 & 1.00 & 1.00 & 1.00 &1.00\\
60&4000&1.00 & 1.00 & 1.00 & 1.00 & 1.00 &1.00\\
100&4000&1.00 & 1.00 & 1.00 & 1.00 & 1.00 &1.00\\
60&8000&1.00 & 1.00 & 1.00 & 1.00 & 1.00 &1.00\\
100&8000&1.00 & 1.00 & 1.00 & 1.00 & 1.00 &1.00\\
\hline
10 & 50 & 8.00 & 8.00 & 8.00& 8.00 & 8.00 & 8.00\\
10 & 100 &8.00 & 8.00 & 8.00& 8.00 & 8.00 & 8.00\\
20 & 100 & 1.01(0.15) & 1.01(0.12) & 1.08(0.53) & 6.96(0.88) & 6.35(0.98) & 7.84(0.40)\\
100 & 10 & 1.08(0.73) & 1.03(0.50) & 1.15(1.01) & 8.00 & 8.00 & 8.00\\
100 & 20 & 1.00(0.03) & 1.00(0.03) & 1.00(0.03) & 5.88(0.76) & 5.12(0.77) & 7.35(0.63)\\
\hline
\end{tabular}
\end{table}

\begin{table}
\caption{Comparison between $PC^\ast_{p1}, PC^\ast_{p2}, PC^\ast_{p3}$ and $PC_{p1}, PC_{p2}, PC_{p3}$ in terms of  the mean estimation numbers of PCs for $m=3, \theta=3$.}
\centering
\tiny
\begin{tabular}{cc|ccc|ccc}
\hline
 N& T & $PC^\ast_{p1}$ & $PC^\ast_{p2}$ & $PC^\ast_{p3}$ & $PC_{p1}$ & $PC_{p2}$ & $PC_{p3}$\\
\hline
100 & 40 & 2.98(0.15) & 2.95(0.22) & 3.00(0.06) &3.00 &3.00&3.90\\
100 & 60 & 3.00(0.03) &3.00(0.04) &3.00& 3.01(0.08) & 3.00 & 4.37(0.64)\\
200 & 60 & 3.00 &3.00 &3.00& 3.00 & 3.00 & 4.18(0.63)\\
500 & 60 & 3.00 &3.00 &3.00& 3.00 & 3.00 & 3.00\\
1000&60&3.00 &3.00 &3.00& 3.00 & 3.00 & 3.00\\
2000 & 60 &3.00 &3.00 &3.00& 3.00 & 3.00 & 3.00\\
100 & 100 &3.00 &3.00 &3.00 & 3.00 & 3.00 & 5.62(0.72)\\
200 & 100 &3.00 &3.00 &3.00& 3.00 & 3.00 & 3.00\\
500 &100&3.00 &3.00 &3.00& 3.00 & 3.00 & 3.00\\
1000&100&3.00 &3.00 &3.00& 3.00 & 3.00 & 3.00\\
2000&60&3.00 &3.00 &3.00& 3.00 & 3.00 & 3.00\\
40 & 100 &2.99(0.10) &2.98(0.14) &3.00 & 3.07(0.26) & 3.01(0.07) & 5.04(0.72) \\
60 & 100 & 3.00 & 3.00(0.03) &3.00 &3.00 & 3.00 & 4.65(0.69)\\
60 & 200 & 3.00 &3.00 &3.00& 3.00 & 3.00 & 3.00\\
60 & 500 & 3.00 &3.00 &3.00& 3.00 & 3.00 & 3.00\\
60 & 1000&3.00 &3.00 &3.00& 3.00 & 3.00 & 3.00\\
60 & 2000&3.00 &3.00 &3.00& 3.00 & 3.00 & 3.00\\
4000&60&3.00 &3.00 &3.00& 3.00 & 3.00 & 3.00\\
4000&100&3.00 &3.00 &3.00& 3.00 & 3.00 & 3.00\\
8000&60&3.00 &3.00 &3.00& 3.00 & 3.00 & 3.00\\
8000&100&3.00 &3.00 &3.00& 3.00 & 3.00 & 3.00\\
60&4000&3.00 &3.00 &3.00& 3.00 & 3.00 & 3.00\\
100&4000&3.00 &3.00 &3.00& 3.00 & 3.00 & 3.00\\
60&8000&3.00 &3.00 &3.00& 3.00 & 3.00 & 3.00\\
100&8000&3.00 &3.00 &3.00& 3.00 & 3.00 & 3.00\\
\hline
10 & 50 & 8.00 & 8.00 & 8.00 & 8.00 & 8.00 & 8.00\\
10 & 100&8.00 & 8.00 & 8.00 & 8.00 & 8.00 & 8.00\\
20 & 100 & 2.89(0.32) & 2.85(0.37) &2.95(0.27) & 6.55(0.74) & 5.96(0.77) & 7.62(0.55)\\
100 & 10 & 2.57(1.35) & 2.43(1.19) &2.77(1.54) & 8.00 & 8.00 & 8.00\\
100 & 20 & 2.46(0.63) & 2.37(0.65) & 2.65(0.52) & 6.15(0.69) & 5.46(0.68) & 7.49(0.59)\\
\hline
\end{tabular}
\end{table}

\begin{table}
\caption{Comparison between $PC^\ast_{p1}, PC^\ast_{p2}, PC^\ast_{p3}$ and $PC_{p1}, PC_{p2}, PC_{p3}$ in terms of  the mean estimation numbers of PCs for $m=5, \theta=5$.}
\centering
\tiny
\begin{tabular}{cc|ccc|ccc}
\hline
 N& T & $PC^\ast_{p1}$ & $PC^\ast_{p2}$ & $PC^\ast_{p3}$ & $PC_{p1}$ & $PC_{p2}$ & $PC_{p3}$\\
\hline
100 & 40 & 3.83(0.77) & 3.49(0.77) & 4.51(0.58) &5.00(0.07) &4.98(0.15)&5.36(0.51)\\
100 & 60 & 4.66(0.50) &4.36(0.61) &4.98(0.13)& 5.00(0.03) & 5.00(0.06) & 5.27(0.45)\\
200 & 60 & 4.95(0.22) &4.90(0.30) &4.99(0.08)& 5.00 & 5.00 & 5.00\\
500 & 60 & 5.00(0.04) &5.00(0.07) &5.00(0.03)& 5.00 & 5.00 & 5.00\\
1000&60&5.00(0.04) &5.00(0.04) &5.00&5.00 & 5.00 & 5.00\\
2000 & 60 &5.00(0.03) &5.00(0.03) &5.00(0.03)& 5.00 & 5.00 & 5.00\\
100 & 100 &4.(0.12) &4.90(0.30) &5.00 & 5.00 & 5.00 & 6.18(0.63)\\
200 & 100 &5.00 &5.00 &5.00& 5.00 & 5.00 & 5.00\\
500 &100&5.00 &5.00 &5.00& 5.00 & 5.00 & 5.00\\
1000&100&5.00 &5.00 &5.00& 5.00 & 5.00 & 5.00\\
2000&60&5.00 &5.00 &5.00& 5.00 & 5.00 & 5.00\\
40 & 100 &4.25(0.68) &3.92(0.75) &4.77(0.44) & 4.98(0.04) & 5.66(0.14) & 5.66(0.57) \\
60 & 100 & 4.76(0.44) & 4.47(0.60) &4.76(0.10) &5.00(0.03) & 4.99(0.08) & 5.46(0.56)\\
60 & 200 & 4.97(0.17) &4.94(0.24) &5.00& 5.00 & 5.00 & 5.00\\
60 & 500 & 5.00(0.05) &5.00(0.06) &5.00(0.04)& 5.00 & 5.00 & 5.00\\
60 & 1000&5.00(0.03) &5.00(0.03) &5.00& 5.00 & 5.00 & 5.00\\
60 & 2000&5.00 &5.00 &5.00& 5.00 & 5.00 & 5.00\\
4000&60&5.00 &5.00 &5.00& 5.00 & 5.00 & 5.00\\
4000&100&5.00 &5.00 &5.00& 5.00 & 5.00 & 5.00\\
8000&60&5.00 &5.00 &5.00& 5.00 & 5.00 & 5.00\\
8000&100&5.00 &5.00 &5.00& 5.00 & 5.00 & 5.00\\
60&4000&5.00 &5.00 &5.00& 5.00 & 5.00 & 5.000\\
100&4000&5.00 &5.00 &5.00& 5.00 & 5.00 & 5.00\\
60&8000&5.00 &5.00 &5.00& 5.00 & 5.00 & 5.00\\
100&8000&5.00 &5.00 &5.00& 5.00 & 5.00 & 5.00\\
\hline
10 & 50 & 8.00 & 8.00 & 8.00 & 8.00 & 8.00 & 8.00\\
10 & 100&8.00 & 8.00 & 8.00 & 8.00 & 8.00 & 8.00\\
20 & 100 & 3.64(0.91) & 3.38(0.94) &4.08(0.79) & 6.65(0.64) & 6.12(0.64) & 7.63(0.51)\\
100 & 10 & 3.10(2.01) & 2.83(1.86) &3.53(2.27) & 8.00 & 8.00 & 8.00\\
100 & 20 & 2.18(0.92) & 1.93(0.92) & 2.65(0.0.90) & 6.56(0.62) & 5.97(0.62) & 7.66(0.50)\\
\hline
\end{tabular}
\end{table}

\begin{table}
\caption{Comparison between $PC^\ast_{p1}, PC^\ast_{p2}, PC^\ast_{p3}$ and $PC_{p1}, PC_{p2}, PC_{p3}$ in terms of  the mean estimation numbers of PCs for $m=5, \theta=3$.}
\centering
\tiny
\begin{tabular}{cc|ccc|ccc}
\hline
 N& T & $PC^\ast_{p1}$ & $PC^\ast_{p2}$ & $PC^\ast_{p3}$ & $PC_{p1}$ & $PC_{p2}$ & $PC_{p3}$\\
\hline
100 & 40 & 4.91(0.30) & 4.81(0.41) & 4.99(0.11) & 5.00(0.03) & 5.00 & 5.59(0.57) \\
100 & 60 & 5.00(0.04) & 4.99(0.11) & 5.00 & 5.00 & 5.00 &5.58(0.57) \\
200 & 60 & 5.00 & 5.00 & 5.00 & 5.00 & 5.00 & 5.00\\
500 & 60 & 5.00 & 5.00 & 5.00 & 5.00 & 5.00 & 5.00\\
1000&60 &5.00 & 5.00 & 5.00 & 5.00 & 5.00 & 5.00\\
2000& 60 & 5.00 & 5.00 & 5.00 & 5.00 & 5.00 & 5.00\\
100 & 100 & 5.00 & 5.00 & 5.00 & 5.00 & 5.00 & 6.84(0.65)\\
200 & 100 & 5.00 & 5.00 & 5.00 & 5.00 & 5.00 & 5.00\\
500 & 100 & 5.00 & 5.00 & 5.00 & 5.00 & 5.00 & 5.00\\
1000&100& 5.00 & 5.00 & 5.00 & 5.00 & 5.00 & 5.00\\
2000&100 & 5.00 & 5.00 & 5.00 & 5.00 & 5.00 & 5.00\\
40 & 100 & 4.97(0.17) & 4.92(0.27) & 5.00(0.04) & 5.02(0.12) & 5.00 & 6.22(0.66)\\
60 & 100 & 5.00(0.04) & 4.99(0.08) & 5.00 & 5.00 & 5.00 & 6.03(0.64)\\
60 & 200 &5.00 & 5.00 & 5.00 & 5.00 & 5.00 & 6.03(0.03)\\
60 & 500 & 5.00 & 5.00 & 5.00 & 5.00 & 5.00 & 5.00\\
60 & 1000 & 5.00 & 5.00 & 5.00 &  5.00 & 5.00 & 5.00\\
60 & 2000 &5.00 & 5.00 & 5.00 &  5.00 & 5.00 & 5.00\\
4000&60 & 5.00 & 5.00 & 5.00 &  5.00 & 5.00 & 5.00\\
4000&100& 5.00 & 5.00 & 5.00 &  5.00 & 5.00 & 5.00\\
8000&60&5.00 & 5.00 & 5.00 &  5.00 & 5.00 & 5.00\\
8000&100&5.00 & 5.00 & 5.00 &  5.00 & 5.00 & 5.00\\
60 & 4000&5.00 & 5.00 & 5.00 &  5.00 & 5.00 & 5.00\\
100 & 4000&5.00 & 5.00 & 5.00 &  5.00 & 5.00 & 5.00\\
60 & 8000&5.00 & 5.00 & 5.00 &  5.00 & 5.00 & 5.00\\
100&8000&5.00 & 5.00 & 5.00 &  5.00 & 5.00 & 5.00\\
\hline
10 & 50 & 8.00 & 8.00 & 8.00 & 8.00 & 8.00 & 8.00\\
10 & 100 &8.00 & 8.00 & 8.00 & 8.00 & 8.00 & 8.00\\
20 & 100 & 4.74(0.51) & 4.62(0.57) & 4.92(0.45) & 7.11(0.63) & 6.65(0.64) & 7.85(0.37)\\
100 & 10 & 4.59(1.99) & 4.35(1.91) & 4.88(2.09) &8.00 & 8.00 & 8.00 \\
100 & 20 & 3.86(0.79) & 3.69(0.81) & 4.13(0.73) & 6.74(0.63) & 6.19(0.62) & 7.77(0.43)\\
\hline
\end{tabular}
\end{table}

\section{Application to the goodness-of-fit test of a PPCA model}

As a third application of the bias-corrected estimator $\widehat{\sigma}_\ast^2$,
we consider the following goodness-of-fit test for the PPCA model (1).
The null hypothesis is then
\[\mathcal{H}_0:~\bSi=\bLa\bLa'+\sigma^2\I_p\text{,}\]
where the number of PCs $m$ is specified.
Following \cite{Anderson1}, the likelihood ratio test (LRT) statistic is
\[T_n = -n L^*\text{,}\]
where
\[L^* = \sum_{j=m+1}^p \log{\frac{\lambda_{n,j}}{\widehat{\sigma}^2}}\text{,}\]
and $\widehat{\sigma}^2$ is the m.l.e.(4) of the variance.
Keeping $p$ fixed while letting $n\rightarrow \infty$, the classical
low-dimensional  theory states that $T_n$ converges to $\chi_q^2$, where $q=p(p+1)/2+m(m-1)/2-pm-1$, see \cite{Anderson1}. However, this classical approximation is again useless  in the large-dimensional
setting. Indeed, it will be shown below  that this criterion leads to
a  high
 false-positive rate. In particular, the test becomes biased since the
size will be much higher than the nominal level
 (see Table \ref{comparison}).

In a way similar to Section 2, we now construct a
corrected version of $T_n$ using Proposition \ref{LRT} and calculus
done in \cite{Bai2} and \cite{Zheng}. As we consider the logarithm of
the eigenvalues of the sample covariance matrix, we will assume in the
sequel that  $p<n$ and $c<1$ to avoid null eigenvalues.

\begin{theorem2} \label{t2}
  Assume the same conditions as in  Theorem 1  and  in
  addition  $c<1$.
  Then, we have
  \[v(c)^{-\frac{1}{2}} \left\{
  L^* - m(c)- ph(c_n)-\eta+(p-m)\log(\beta) \right\}
  ~\overset{\mathcal{D}}{\longrightarrow} ~ \mathcal{N}(0,1)\text{,}\]
where
\begin{eqnarray*}
  m(c) & = & \frac{\log{(1-c)}}{2}~, \qquad
  h(c_n)=\frac{c_n-1}{c_n} \log(1-c_n)-1 ~,\\
  \eta&=&\sum_{i=1}^m \log (1+c\sigma^2\alpha_i^{-1}) ~, \qquad
  \beta=1-\frac{c}{p-m}(m+\sigma^2\sum_{i=1}^m \alpha_i^{-1})~,\\
  v(c)&=&-2\log(1-c)+\frac{2c}{\beta}\left (\frac{1}{\beta}-2 \right )~.
\end{eqnarray*}
\end{theorem2}

Note that the above statistic depends on the unknown  variance
$\sigma^2$ and  the spike eigenvalues $(\alpha_i)$.
First of all, as explained  in Section 2,
consistent estimates of $(\alpha_i)$ are available. By using these
estimates and substituting bias-corrected estimate
$\widehat\sigma_*^2$ for $\sigma^2$, we obtain
consistent estimates
$\widehat v(c_n)$, $\widehat \eta$ and  $\widehat \beta$
of   $ v(c)$, $\eta$ and  $\beta$, respectively. Therefore,
to test $\mathcal{H}_0$,  it is natural to use the statistic
\[ \Delta_n:=
  \widehat v(c_n)^{-\frac{1}{2}}(L^*-m(c_n)-ph(c_n)-
  \widehat \eta+(p-m)\log(\widehat \beta)) ~\text{.} \label{statclrt}
\]
Since $\Delta_n$ is asymptotically standard normal,
the critical region $\{\Delta_n> q_\alpha \}$
where $q_\alpha$ is the $\alpha$th upper quantile of the standard normal, will have an asymptotic  size $\alpha$.
This test will be hereafter referred as the corrected likelihood ratio test (CLRT in short).

\subsection{Monte-Carlo experiments}
We consider again Models 1 and 2 described in Section 2, and a new one (Model~4):
\begin{list}{$\bullet$}{\leftmargin=2em}
\item Model 1: $\mbox{spec}(\bSi)= (25,16,9,0,\dots,0)+\sigma^2(1,\dots,1)$, $\sigma^2=4$, $c=0.9$;
\item Model 2: $\mbox{spec}(\bSi)= (4,3,0,\dots,0)+\sigma^2(1,\dots,1)$, $\sigma^2=2$, $c=0.2$;
\item Model 4: $\mbox{spec}(\bSi)= (8,7,0,\dots,0)+\sigma^2(1,\dots,1)$, $\sigma^2=1$, varying $c$.
\end{list}

Table \ref{t3c6} gives the empirical sizes of the classical
likelihood ratio test (LRT) and the new corrected likelihood ratio test
(CLRT) above. For the LRT, we use the correction proposed by
\cite{Bartlett}, that is replacing $T_n=-nL^*$ by
$\tilde{T}_n=-(n-(2p+11)/6-2m/3)L^*$. The
computations are done under 10000 independent replications and the nominal test level is $0.05$.

{
\begin{table}[!ht]
\begin{center}
\small
\caption{{  Comparison of the empirical size of the classical likelihood ratio test (LRT) and the corrected likelihood ratio test (CLRT) in various settings.}} \label{t3c6}
{ \begin{tabular}{|cll|cc|}
\hline
\multicolumn{3}{|c|}{Settings} & Empirical size of CLRT & Empirical size of LRT \\
\hline
\multirow{3}{*}{Model 1} & $p=90$ & $n=100$ & 0.0497 & 0.9995\\
 & $p=180$ & $n=200$ & 0.0491 & 1\\
 & $p=720$ & $n=800$ & 0.0496 & 1\\
\hline
\multirow{3}{*}{Model 2} & $p=20$ & $n=100$ & 0.0324 & 0.0294\\
 & $p=80$ & $n=400$ & 0.0507 & 0.0390\\
 & $p=200$ & $n=1000$ & 0.0541 & 0.0552\\
\hline
\multirow{4}{*}{Model 4} & $p=5$ & $n=500$ & 0.0108 & 0.0483\\
 & $p=10$ & $n=500$ & 0.0190 & 0.0465\\
 & $p=50$ & $n=500$ & 0.0424 & 0.0445\\
 & $p=100$ & $n=500$ & 0.0459 & 0.0461\\
 & $p=200$ & $n=500$ & 0.0491 & 0.2212\\
 & $p=250$ & $n=500$ & 0.0492 & 0.7395\\
 & $p=300$ & $n=500$ & 0.0509 & 0.9994\\
\hline
\end{tabular} \label{comparison}}
\end{center}
\end{table}}

The empirical sizes of the new CLRT are very close to the nominal one, except
when the ratio $p/n$ is very small (less than 0.1). On the contrary, the empirical sizes of the classical LRT are much higher than the nominal level especially when $c$ is not too small, and the test will always reject the null hypothesis when $p$ becomes large. In particular when $p/n\ge \frac12$, the LRT test tends to reject automatically the null.

\section{Proofs}
Before giving the proofs, we first recall some important results from the random matrix theory which laid the foundation for the proofs of the main results of the paper.

\subsection{Useful results from random matrix theory}
Random matrix theory has become a powerful tool to address new
inference problems in high-dimensional data. For general
background and references, we refer to review papers
\citet{JohnstoneICM} and \citet{JohnTitt09}.

Let $H$ be a probability measure on $\R^+$ and $c>0$ a constant. We define the map
\begin{eqnarray}
g(s)=g_{c,H}(s)= \frac{1}{s} + c \int \frac{t}{1+ts} \, \dif H(t) \label{fund_eq}
\end{eqnarray}
in the set $\mathbb{C}^+ = \{ z \in \mathbb{C}:\Im z>0\}$. The map $g$ is a one-to-one mapping from $\mathbb{C}^+$ onto itself (see \cite{Bai4}, Chapter 6), and the inverse map $m=g^{-1}$ satisfies all the requirements of the Stieltjes transform of a probability measure on $[0,\infty)$. We call this measure $\emph{\b{F}}_{c,H}$. Next, a companion measure $F_{c,H}$ is introduced by the equation $c F_{c,H}=\left ( c-1\right ) \delta_0 + \emph{\b{F}}_{c,H}$ (note that in this equation,  measures can be signed). The measure $F_{c,H}$ is referred as the generalized Mar{\v c}enko-Pastur distribution with index $(c,H)$.

Let $F_n= \frac{1}{p}\sum_{i=1}^p \delta_{ \lambda_{n,i}}$ be the
empirical spectral distribution (ESD) of the sample covariance matrix $\mathbf{S}_n$ defined in (3) with the $\{\lambda_{n,i}\}$ denoting its eigenvalues. Then, it is well-known that under suitable moment conditions, $F_n$ converges to the Mar{\v c}enko-Pastur distribution of index
$(c,\delta_{\sigma^2})$,
simply denoted as $F_{c,\sigma^2}$, with the following density function
\[p_{c,\sigma^2}(x)=
\begin{cases}
  \frac{1}{2\pi x c \sigma^2} \sqrt{\{b(c)-x\}\{x-a(c)\}}~,
  &  \quad  a(c)\le x\le b(c) ~,\\
  0~, &  \quad \mbox{ otherwise.}
\end{cases}
\]
The distribution has an additional mass $(1-1/c)$ at the origin if
$c>1$.

The ESD $H_n$ of $\bSi$ is
\[H_n=\frac{p-m}{p} \delta_{\sigma^2} + \frac{1}{p} \sum_{i=1}^m \delta_{\alpha_i+\sigma^2}~\text{,}\]
and $H_n \rightarrow \delta_{\sigma^2}$. Define the normalized empirical process
\[G_n(f)= p \int_{\R} f(x)[F_n-F_{c_n,H_n}] (\dif x)\text{, } f \in \mathcal{A}\text{,}\]
where $\mathcal{A}$ is the set of analytic functions $f:\mathcal{U}
\to \mathbb{C}$, with $\mathcal{U}$ an open set of $\mathbb{C}$ such
that $[\text{\large 1}_{(0,1)}(c)a(c),b(c)] \subset \mathcal{U}$.
We will  need the
following CLT  which is a combination of
Theorem 1.1 of \citet{BS04} and a recent addition  proposed in
\citet{zheng2014substitution}.

\begin{proposition} \label{LRT}
 We assume the same conditions as in Theorem 1.
 Then, for any functions
 $f_1,\dots,f_k \in \mathcal{A}$,
 the random
 vector $(G_n(f_1),\dots,G_n(f_k))$ converges to a $k$-dimensional
 Gaussian vector with mean
 vector
\[m(f_j) = \frac{f_j(a(c))+f_j(b(c))}{4}- \frac{1}{2\pi} \int_{a(c)}^{b(c)} \frac{f_j(x)}{\sqrt{4c\sigma^4-(x-\sigma^2-c\sigma^2)^2}}\, \dif x\text{, }j=1,\dots,k\text{,}\]
and covariance function
\begin{eqnarray}
v(f_j,f_l)&=&-\frac{1}{2\pi^2}\oint_{\mathcal{C}_1} \oint_{\mathcal{C}_2} \frac{f_j(z_1)f_l(z_2)}{(\underline{m}(z_1)-\underline{m}(z_2))^2}\, \dif \underline{m}(z_1)\dif \underline{m}(z_2)\text{, }j,l=1,\dots,k\text{,} \label{cov}
\end{eqnarray}
where $\underline{m}(z)$ is the Stieltjes transform of
$\b{F}_{c,{\sigma^2}}= (1-c) \delta_0+cF_{c,{\sigma^2}}$.
The contours $\mathcal{C}_1$ and $\mathcal{C}_2$ are non overlapping and both contain the support of $F_{c,{\sigma^2}}$.
\end{proposition}

An important and subtle point here is that the centering term in
$G_n(f)$ in the above CLT is defined with respect to the
Marc\v{c}enko-Pastur distribution   $F_{c_n,H_n}$ with ``current''
index $(c_n,H_n)$ instead of the limiting distribution
$F_{c,{\sigma^2}}$ with
index $(c,{\sigma^2})$.
In contrast, the limiting mean function $m(f_j)$ and
covariance function $v(f_j,f_l)$ depend on the limiting distribution $F_{c,{\sigma^2}}$
only.

\subsection{Proof of Theorem 1}
We have
  \begin{eqnarray*}
    (p-m) \widehat{\sigma}^2 &=& 		\sum_{i=1}^p \lambda_{n,i} - \sum_{i=1}^m \lambda_{n,i}\text{.}
  \end{eqnarray*}
  By (\ref{llambda}),
  \begin{eqnarray}
    \sum_{i=1}^m \lambda_{n,i} &\longrightarrow& \sum_{i=1}^m \left (  \alpha_i + \frac{c\sigma^4}{\alpha_i}\right ) +\sigma^2 m(1+c) \text{ a.s.}\label{1}
  \end{eqnarray}
  For the first term, we have
  \begin{eqnarray*}
    \sum_{i=1}^p \lambda_i &=& p \int x d F_n(x)\\
	&=& p \int x \, \dif (F_n-F_{c_n,H_n})(x) + p \int x \, \dif F_{c_n,H_n}(x)\\
	&=& G_n(x) + p \int x \, \dif F_{c_n,H_n}(x) \text{.}
  \end{eqnarray*}
  By Proposition \ref{LRT}, the first term is asymptotically normal
  \begin{eqnarray*}
    G_n(x)= \sum_{i=1}^p \lambda_{n,i} - p \int x  \, \dif F_{c_n,H_n}(x) \overset{\mathcal{D}}{\longrightarrow} \mathcal{N}(m(x),v(x))\text{,}
  \end{eqnarray*}
  with asymptotic  mean
  \begin{eqnarray}
    m(x)&=&0~, \label{mx}
  \end{eqnarray}
and asymptotic variance
\begin{eqnarray}
  v(x)&=&2c  \sigma^4 ~.\label{vx}
\end{eqnarray}
The derivation of these two formula are given in the Section 6.
Furthermore, by Lemma 1 of \cite{Bai3},
\begin{eqnarray*}
\int x \, \dif F_{c_n,H_n}(x) &=& \int t \, \dif H_n(t)
= \sigma^2 + \frac{1}{p} \sum_{i=1}^m\alpha_i\text{.}
\end{eqnarray*}
So we have
\begin{eqnarray}
\sum_{i=1}^p \lambda_{n,i} - p\sigma^2 - \sum_{i=1}^m \alpha_i &\overset{\mathcal{D}}{\longrightarrow}& \mathcal{N}(0,2c\sigma^4)\text{.}\label{2}
\end{eqnarray}
By (\ref{1}) and (\ref{2}) and using Slutsky's lemma, we obtain
\[(p-m)(\widehat{\sigma}^2-\sigma^2) + c\sigma^2 \left ( m + \sigma^2 \sum_{i=1}^m \frac{1}{\alpha_i}\right ) \overset{\mathcal{D}}{\longrightarrow} \mathcal{N}(0,2c\sigma^4)\text{.}\]

\subsection{Proof of Theorem 2}
We have
\begin{eqnarray*}
\frac{p-m}{\sigma^2\sqrt{2c_n}}\left(\widehat{\sigma}_\ast^2 -\sigma^2\right) &=& \frac{p-m}{\sigma^2\sqrt{2c_n}}\left(\widehat{\sigma}^2- \sigma^2\right)+b\left(\widehat{\sigma}^2\right) \frac{\widehat{\sigma}^2}{\sigma^2}\\
&=& \left\{\frac{p-m}{\sigma^2\sqrt{2c_n}}\left(\widehat{\sigma}^2- \sigma^2\right)+b(\sigma^2)\right\}+\frac{1}{\sigma^2} \left\{b\left(\widehat{\sigma}^2\right)\widehat{\sigma}^2 -b(\sigma^2)\sigma^2\right\}.
\end{eqnarray*}
Since $ \widehat{\sigma}^2\overset{\mathcal{P}}{\longrightarrow} \sigma^2$, by continuity, the second expression tends to 0 in probability and the conclusion follows from Theorem 1.

\subsection{Proof of Theorem 3}

We have
\begin{eqnarray*}
L^* &=& \sum_{i=m+1}^p \log \frac{\lambda_{n,i}}{\widehat{\sigma}^2}\\
    &=& \sum_{i=m+1}^p \log \frac{\lambda_{n,i}}{\sigma^2} - \sum_{i=m+1}^p \log \frac{\widehat{\sigma}^2}{\sigma^2}\\
    &=& \sum_{i=m+1}^p \log \frac{\lambda_{n,i}}{\sigma^2} - (p-m)\log \left ( \frac{1}{p-m}\sum_{i=m+1}^p \frac{\lambda_{n,i}}{\sigma^2} \right ) \\
    &=& L_1-(p-m)\log \left( \frac{L_2}{p-m}\right )\text{,}
\end{eqnarray*}
where we have defined a two-dimensional vector $(L_1,L_2)=(\sum_{i=m+1}^p \log \frac{\lambda_{n,i}}{\sigma^2},\sum_{i=m+1}^p  \frac{\lambda_{n,i}}{\sigma^2})$.

\vspace{4mm}
\noindent {\bf CLT when $\sigma^2=1$.} To start with, we consider the case $\sigma^2=1$.
We have
\begin{eqnarray*}
  L_1&=& p \int \log(x) \, \dif F_n(x) - \sum_{i=1}^m \log{\lambda_{n,i}}\\
  &=& p \int \log(x) \, \dif (F_n-F_{c_n,H_n})(x) + p \int \log(x) \, \dif F_{c_n,H_n}(x)- \sum_{i=1}^m \log \lambda_{n,i}\text{.}
\end{eqnarray*}
Similarly, we have
\[L_2= p \int x \, \dif (F_n-F_{c_n,H_n})(x) + p \int x \, \dif F_{c_n,H_n}(x)- \sum_{i=1}^m \lambda_{n,i}\text{.}\]
By Proposition \ref{LRT}, we find that
\begin{eqnarray}
  p \left ( \begin{array}{l}
    \int \log(x) \, \dif (F_n-F_{c_n,H_n})(x)\\
	\int x \, \dif (F_n-F_{c_n,H_n})(x) \end{array} \right )
  &\overset{\mathcal{D}}{\longrightarrow}& \mathcal{N} \left ( \left (\begin{array}{c}
    m_1(c)\\
	m_2(c) \end{array} \right )
  , \left (  \begin{array}{lr} v_1(c) & v_{1,2} (c)\\
    v_{1,2}(c) & v_2(c)\end{array}\right ) \right ) \label{CLTBS2}
\end{eqnarray}
with $m_2(c)=0$ and $v_2(c)=2c$ and
\begin{eqnarray}
m_1(c)&=&\frac{\log{(1-c)}}{2}\text{,}\label{mlogx}
\end{eqnarray}
\begin{eqnarray}
v_1(c)&=&-2\log{(1-c)}\text{,}\label{vlogx}
\end{eqnarray}
\begin{eqnarray}
v_{1,2}(c)&=&2c\text{.}\label{v12}
\end{eqnarray}
Formulae of $m_2$ and $v_2$ have been established in the proof of
Theorem 1 and
the others are derived  in next subsection.

In Theorem \ref{T1}, with $\sigma^2=1$, we found that
\[\int x \, \dif F_{c_n,H_n}(x)=1 + \frac{1}{p} \sum_{i=1}^m \alpha_i\text{,}\]
and
\[\sum_{i=1}^m \lambda_{n,i} \overset{\mbox{a.s.}}\longrightarrow \sum_{i=1}^m \left (  \alpha_i + \frac{c}{\alpha_i}\right ) + m(1+c) \text{.}\]
For the last term of $L_1$, by (\ref{llambda}), we have
\[ \log \lambda_{n,i} \longrightarrow  \log (\phi(\alpha_i+1))=\log \left (  (\alpha_i + 1)  ( 1+ c\alpha_i^{-1})\right ) \text{ a.s.}\]
Furthermore, by \citet{Wang}, we have
\[\int \log(x) \, \dif F_{c_n,H_n}(x)=\frac{1}{p} \sum_{i=1}^m \log(\alpha_i+1) +h(c_n) + o\left ( \frac{1}{p}\right )\text{,}\]
where
\begin{eqnarray}
h(c_n)&=&\int \log(x) dF_{c_n,\delta_{1}}(x)
=\frac{c_n-1}{c_n} \log(1-c_n)-1\text{.}\label{hc}
\end{eqnarray}
can be calculated using the density of the Mar{\v c}enko-Pastur law (see 6.4).
Summarising, we have obtained that
\[L_1-m_1(c)-ph(c_n)+\eta(c,\alpha)\overset{\mathcal{D}}{\longrightarrow} \mathcal{N} \left ( 0,v_1(c) \right )\text{,}\]
where $h(c_n)=\frac{c_n-1}{c_n} \log(1-c_n)-1$ and $\eta(c,\alpha)=\sum_{i=1}^m \log (1+c\sigma^2\alpha_i^{-1})$. Similarly, we have
\[ L_2 - (p-m) + \rho(c,\alpha)\overset{\mathcal{D}}{\longrightarrow} \mathcal{N} \left ( 0,v_2(c) \right )\text{,}\]
where  $\rho(c,\alpha)=c(m+\sum_{i=1}^m \alpha_i^{-1})$.

Using \eqref{CLTBS2} and Slutsky's lemma,
\[\left ( \begin{array}{l}
                  L_1\\
				  L_2 \end{array} \right )
\overset{\mathcal{D}}{\longrightarrow} \mathcal{N} \left ( \left (\begin{array}{c}
                                                               m_1(c)+ph(c_n)-\eta(c,\alpha)\\
                                     p-m-\rho(c,\alpha) \end{array} \right )
, \left (  \begin{array}{lr} v_1(c) & v_{1,2} (c)\\
v_{1,2}(c_n) & v_2(c_n)\end{array}\right ) \right )\text{,}\]
with $h(c_n)=\frac{c_n-1}{c_n} \log(1-c_n)-1$, $\eta(c,\alpha)=\sum_{i=1}^m \log (1+c\sigma^2\alpha_i^{-1})$ and $\rho(c,\alpha)=c(m+\sum_{i=1}^m \alpha_i^{-1})$.
\vspace{4mm}

\noindent{\bf CLT with general $\sigma^2$.} When $\sigma^2=1$,
\[\mbox{spec}(\bSi)=(\alpha_1+1,\dots,\alpha_m+1,1,\dots,1)\text{,}\]
whereas in the general case
\begin{eqnarray*}
\mbox{spec}(\bSi)&=&(\alpha_1+\sigma^2,\dots,\alpha_m+\sigma^2,\sigma^2,\dots,\sigma^2)\\
				   &=&\sigma^2 \left (\frac{\alpha_1}{\sigma^2}+1,\dots,\frac{\alpha_m}{\sigma^2}+1,\dots,1 \right )\text{.}
\end{eqnarray*}
Thus, if we consider $\lambda_i/\sigma^2$, we will find the same CLT by replacing the $(\alpha_i)_{1 \le i \le m}$ by $\alpha_i/\sigma^2$. Furthermore, we divide $L_2$ by $p-m$ to find
\begin{eqnarray*} \left ( \begin{array}{c}
                  L_1\\
				  \frac{L_2}{p-m} \end{array} \right )
&\overset{\mathcal{D}}{\longrightarrow}& \mathcal{N} \left ( \left (\begin{array}{c}
                                                               m_1(c)+ph(c_n)-\eta(c,\alpha/\sigma^2)\\
                                     1-\frac{\rho(c,\alpha/\sigma^2)}{p-m} \end{array} \right )
, \left (  \begin{array}{cc} \frac{2c}{(p-m)^2} & \frac{2c}{p-m}\\
\frac{2c}{p-m} & -2\log(1-c)\end{array}\right ) \right )\text{,}\label{CLT_CLRT}
\end{eqnarray*}
with $\eta(c,\alpha/\sigma^2)=\sum_{i=1}^m \log (1+c\sigma^2\alpha_i^{-1})$, $\rho(c,\alpha/\sigma^2)=c(m+\sigma^2\sum_{i=1}^m \alpha_i^{-1})$ and $h(c_n)=\frac{c_n-1}{c_n} \log(1-c_n)-1$.

\vspace{4mm}
\noindent {\bf Asymptotic distribution of $L^*$}. We have $L^*=g(L_1,L_2/(p-m))$, with $g(x,y)=x-(p-m)\log(y)$. We will apply the multivariate delta-method on (\ref{CLT_CLRT}) with the function $g$. We have $\bigtriangledown g(x,y)=\left ( 1,-\frac{p-m}{y}\right )$ and
\[L^* \overset{\mathcal{D}}{\longrightarrow} \mathcal{N} ( \beta_1-(p-m)\log(\beta_2),\bigtriangledown g(\beta_1,\beta_2) \mbox{ cov}(L_1,L_2/(p-m)) \bigtriangledown g(\beta_1,\beta_2)')\text{,} \]
with $\beta_1=m_1(c)+ph(c_n)-\eta(c,\alpha/\sigma^2)$ and
$\beta_2=1-\frac{\rho(c,\alpha/\sigma^2)}{p-m}$. After some standard
calculation,  we finally find
\[L^* \overset{\mathcal{D}}{\longrightarrow} \mathcal{N} \left
(m_1(c)+ph(c_n)-\eta \left (c,\frac{\alpha}{\sigma^2} \right
)-(p-m)\log(\beta_2),-2\log(1-c)+ \frac{2c}{\beta_2}\left (\frac{1}{\beta_2}-2 \right ) \right )\text{.}\]

\subsection{Complementary proofs}
\subsubsection*{Proof of (\ref{sigmaMLE})}

The general theory of the m.l.e. for the PPCA model~\eqref{model}
in the classical setting has been
developed in \citet{Anderson2} with in particular  the following
result.

\begin{proposition} \label{P1}
Let $\Theta = (\theta_{ij})_{1 \le i,j \le p} = \bPsi -
\bLa(\bLa'\bPsi^{-1}\bLa)^{-1}\bLa'$. If $(\theta_{ij}^2)_{1 \le i,j
  \le p}$ is nonsingular, if $\bLa$ and $\bPsi$ are identified by the
condition that $\bLa'\bPsi \bLa$ is diagonal and the diagonal elements
are different and ordered, if $\mathbf{S}_n \rightarrow \bLa\bLa'+
\bPsi$ in probability and if $\sqrt{n}(\mathbf{S}_n-\bSi)$ has a
limiting distribution, then $\sqrt{n}(\widehat{\bLa} - \bLa)$ and
$\sqrt{n}(\widehat{\bPsi} - \bPsi)$ have a limiting distribution. The
covariance of $\sqrt{n}(\widehat{\bPsi}_{ii} - \bPsi_{ii})$ and
$\sqrt{n}(\widehat{\bPsi}_{jj} - \bPsi_{jj})$ in the limiting
distribution is $2 \bPsi_{ii}^2\bPsi_{jj}^2\xi^{ij}$ $(1 \le i,j \le
p)$,  where $(\xi^{ij})=(\theta_{ij}^2)^{-1}$.
\end{proposition}

To prove the CLT (\ref{sigmaMLE}),
by Proposition \ref{P1}, we know that the inverse of the Fisher
information matrix  is
$\mathcal{I}^{-1}(\psi_{11},\dots,\psi_{pp})=(2\psi_{ii}^2\psi_{jj}^2\xi^{ij})_{ij}$. We
have to change the parametrization: in our case, we have
$\psi_{11}=\dots=\psi_{pp}$. Let $g:\R \rightarrow \R^p$, $a \mapsto
(a,\dots,a)$.  The information matrix  in this new parametrization becomes
\[\mathcal{I}(\sigma^2)=J'\mathcal{I}(g(\sigma^2))J\text{,}\]
where $J$ is the Jacobian matrix of $g$. As
\[\mathcal{I}(g(\sigma^2))=\frac{1}{2\sigma^8}(\theta_{ij}^2)_{ij}\text{,}\]
we  have
\[\mathcal{I}(\sigma^2)=\frac{1}{2\sigma^8} \sum_{i,j=1}^p \theta_{ij}^2\text{,}\]
and
\begin{eqnarray*}
\Theta=(\theta_{ij})_{ij} &=& \bPsi - \bLa(\bLa'\bPsi^{-1}\bLa)^{-1}\bLa'\\
			&=& \sigma^2(\I_p-\bLa(\bLa'\bLa)^{-1}\bLa')\text{.}
\end{eqnarray*}
By hypothesis, we have $\bLa'\bLa=\mbox{diag}(d_1^2,\dots,d_m^2)$. Consider the Singular Value Decomposition of $\bLa$, $\bLa=\U \D \mathbf{V}$, where $\U$ is a $p\times p$ matrix such that $\U\U'=\I_p$, $\mathbf{V}$ is a $m\times m$ matrix such that $\mathbf{V}'\mathbf{V}=\I_m$, and $\D$ is a $p\times m$ diagonal matrix with $d_1,\dots,d_m$ as diagonal elements. As $\bLa'\bLa$ is diagonal, $\mathbf{V}=\I_m$, so $\bLa=\U\D$. By elementary calculus, one can find that
\[\bLa(\bLa'\bLa)^{-1}\bLa'= \mbox{diag}(\underbrace{1,\dots,1}_{m},\underbrace{0,\dots,0}_{p-m})\text{,}\]
so
\[\Theta=\sigma^2\mbox{diag}(\underbrace{0,\dots,0}_{m},\underbrace{1,\dots,1}_{p-m})\text{.}\]
Finally,
\[\mathcal{I}(\sigma^2)=\frac{1}{2\sigma^8}(p-m)\sigma^4=\frac{p-m}{2\sigma^4}~,\]
and the asymptotic variance of $\widehat\sigma^2$ is
\[s^2=\mathcal{I}^{-1}(\sigma^2)=\frac{2\sigma^4}{p-m}\text{.}\]

\subsubsection*{Proof of (\ref{mx})}
By Proposition \ref{LRT}, for $g(x)=x$, by using the variable change $x=\sigma^2(1+c-2\sqrt{c}\cos\theta)$, $0 \le \theta \le \pi$, we have
\begin{eqnarray*}
m(g) &=& \frac{g(a(c))+g(b(c))}{4}- \frac{1}{2\pi} \int_{a(c)}^{b(c)} \frac{x}{\sqrt{4c\sigma^4-(x-\sigma^2-c\sigma^2)^2}} \, \dif x\text{, }j=1,\dots,k \\
&=& \frac{\sigma^2(1+c)}{2}-\frac{\sigma^2}{2\pi} \int_0^{\pi}(1+c-2\sqrt{c}\cos\theta) \, \dif \theta\\
&=&0\text{.}
\end{eqnarray*}

\subsubsection*{Proof of (\ref{vx})}
Let $\underline{s}(z)$ be the Stieltjes transform of
$(1-c)\text{\large 1}_{[0,\infty)} +cF_{c,\delta_1}$. One can show that
\[\underline{m}(z)=\frac{1}{\sigma^2}\underline{s}\left(\frac{z}{\sigma^2}\right)\text{.}\]
Then, in Proposition \ref{LRT}, we have
\begin{eqnarray}
v(f_j,f_l)&=&-\frac{1}{2\pi^2}\oint \oint \frac{f_j(\sigma^2z_1)f_l(\sigma^2z_2)}{(\underline{s}(z_1)-\underline{s}(z_2))^2}\, \dif \underline{s}(z_1)\, \dif \underline{s}(z_2)\text{, }j,l=1,\dots,k\text{.} \label{transfo}
\end{eqnarray}
For $g(x)=x$, we have
\begin{eqnarray*}
v(g)&=&-\frac{1}{2\pi^2}\oint \oint \frac{g(\sigma^2 z_1)g(\sigma^2z_2)}{(\underline{s}(z_1)-\underline{s}(z_2))^2}\, \dif \underline{s}(z_1)\, \dif \underline{s}(z_2)\\
    &=&-\frac{\sigma^4}{2\pi^2}\oint \oint \frac{z_1z_2}{(\underline{s}(z_1)-\underline{s}(z_2))^2}\, \dif \underline{s}(z_1)\, \dif \underline{s}(z_2)\\
    &=&2c\sigma^4\text{,}
\end{eqnarray*}
where $-\frac{1}{2\pi^2}\oint \oint \frac{z_1z_2}{(\underline{s}(z_1)-\underline{s}(z_2))^2} \, \dif \underline{s}(z_1)\, \dif \underline{s}(z_2)=2c$ is calculated in \cite{Bai2} (it corresponds to $v(z_1,z_2)$, Section 5, proof of (3.4)).

\subsubsection*{Proof of (\ref{mlogx})}
By Proposition \ref{LRT}, for $\sigma^2=1$ and $g(x)=\log(x)$, by using the variable change $x=1+c-2\sqrt{c}\cos\theta$, $0 \le \theta \le \pi$, we have
\begin{eqnarray*}
m(g) &=& \frac{g(a(c))+g(b(c))}{4}- \frac{1}{2\pi} \int_{a(c)}^{b(c)} \frac{x}{\sqrt{4c-(x-1-c)^2}} \, \dif x\text{, }j=1,\dots,k \\
&=& \frac{\log(1-c)}{2}-\frac{1}{2\pi} \int_0^{\pi}\log(1+c-2\sqrt{c}\cos\theta)\, \dif \theta\\
&=& \frac{\log(1-c)}{2}-\frac{1}{4\pi} \int_0^{2\pi}\log|1-\sqrt{c}e^{i\theta}|^2 \, \dif \theta\\
&=& \frac{\log(1-c)}{2}\text{,}
\end{eqnarray*}
where $\int_0^{2\pi}\log|1-\sqrt{c}e^{i\theta}|^2 \, \dif \theta=0$ is calculated in \cite{Bai4}.

\subsubsection*{Proof of (\ref{vlogx})}
By Proposition \ref{LRT} and (\ref{transfo}), for $\sigma^2=1$ and $g(x)=x$, we have
\begin{eqnarray*}
v(g)&=&-\frac{1}{2\pi^2}\oint \oint \frac{g(z_1)g(z_2)}{(\underline{s}(z_1)-\underline{s}(z_2))^2} \, \dif \underline{s}(z_1) \, \dif \underline{s}(z_2)\\
    &=&-\frac{1}{2\pi^2}\oint \oint \frac{\log(z_1)\log(z_2)}{(\underline{s}(z_1)-\underline{s}(z_2))^2}\, \dif \underline{s}(z_1)\dif \underline{s}(z_2)\\
    &=&-2\log(1-c_n)\text{,}
\end{eqnarray*}
where the last integral is calculated in \cite{Bai4}.

\subsubsection*{Proof of (\ref{hc})}
$F_{c_n,\delta_1}$ is the Mar{\v c}enko-Pastur distribution of index $c_n$. By using the variable change $x=1+c_n-2\sqrt{c_n}\cos\theta$, $0 \le \theta \le \pi$, we have
\begin{eqnarray*}
\int \log(x) dF_{c_n,\delta_{1}}(x) &=& \int_{a(c_n)}^{b(c_n)}\frac{\log x}{2\pi x c_n} \sqrt{(b(c_n)-x)(x-a(c_n))}\, \dif x\\
&=& \frac{1}{2\pi c_n} \int_0^{\pi} \frac{\log(1+c_n-2\sqrt{c_n}\cos\theta)}{1+c_n-2\sqrt{c_n}\cos\theta}4c_n \sin^2\theta \, \dif \theta\\
&=& \frac{1}{2\pi} \int_0^{2\pi} \frac{2\sin^2\theta}{1+c_n-2\sqrt{c_n}\cos\theta}\log|1-\sqrt{c_n}e^{i\theta}|^2  \, \dif \theta\\
&=& \frac{c_n-1}{c_n} \log(1-c_n)-1\text{,}
\end{eqnarray*}
where the last integral is calculated in \cite{Bai4}.

\subsubsection*{Proof of (\ref{v12})}
In the normal case with $\sigma^2=1$, \cite{Zheng} gives the following equivalent expression of (\ref{cov}):
\[v(f_j,f_l)=-\underset{r \rightarrow 1^+}{\lim} \frac{\kappa}{4\pi^2}\oint\oint_{|\xi_1|=|\xi_2|=1} f_j(|1+h\xi_1|^2)f_l(|1+h\xi_2|^2)\frac{1}{(\xi_1-r\xi_2)^2} \, \dif \xi_1 \, \dif \xi_2\text{,}\]
where $\kappa=2$ in the real case and $h=\sqrt{c}$ in our case. We take $f_j(x)=\log(x)$ and $f_l(x)=x$, so we need to calculate
\[v(\log(x),x)=-\underset{r \rightarrow 1^+}{\lim} \frac{1}{2\pi^2}\oint\oint_{|\xi_1|=|\xi_2|=1} |1+\sqrt{c}\xi_2|^2 \frac{\log(|1+\sqrt{c}\xi_1|^2)}{(\xi_1-r\xi_2)^2} \, \dif \xi_1 \, \dif \xi_2\text{.}\]
We follow the calculations done in \cite{Zheng}: when $|\xi|=1$, $|1+\sqrt{c}\xi|^2=(1+\sqrt{c}\xi)(1+\sqrt{c}\xi^{-1})$, so $\log(|1+\sqrt{c}\xi|^2)=\frac{1}{2}\left ( \log(1+\sqrt{c}\xi)^2+\log(1+\sqrt{c}\xi^{-1})^2 \right )$. Consequently,
\begin{eqnarray*}
\oint_{|\xi_1|=1} \frac{\log(|1+\sqrt{c}\xi_1|^2)}{(\xi_1-r\xi_2)^2} \, \dif \xi_1
&=& \frac{1}{2}\oint_{|\xi_1|=1}\frac{\log(1+\sqrt{c}\xi_1)^2}{(\xi_1-r\xi_2)^2} \, \dif \xi_1 + \frac{1}{2}\oint_{|\xi_1|=1}\frac{\log(1+\sqrt{c}\xi_1^{-1})^2}{(\xi_1-r\xi_2)^2} \, \dif \xi_1\\
&=& \frac{1}{2}\oint_{|\xi_1|=1}\log(1+\sqrt{c}\xi_1)^2 \left ( \frac{1}{(\xi_1-r\xi_2)^2}+\frac{1}{(1-r\xi_1\xi_2)^2} \right )\, \dif \xi_1\\
&=& 0 +i \pi \left ( \frac{1}{(r\xi_2)^2}  \frac{2\sqrt{c}}{1+\frac{\sqrt{c}}{r\xi_2}} \right )\\
&=& 2 i \pi \frac{\sqrt{c}}{r\xi_2 (r\xi_2+\sqrt{c})}\text{.}
\end{eqnarray*}
Thus,
\begin{eqnarray*}
v(\log(x),x)&=&\frac{1}{i\pi} \oint_{|\xi_2|=1} |1+\sqrt{c}\xi_2|^2 \frac{\sqrt{c}}{\xi_2 (\xi_2+\sqrt{c})} \, \dif \xi_2\\
			&=&\frac{1}{i\pi} \oint_{|\xi|=1} \left (1+c+c ( \xi+\xi^{-1} ) \right) \frac{\sqrt{c}}{\xi (\xi+\sqrt{c})} \, \dif \xi\\
			&=&\frac{1}{i\pi} \oint_{|\xi|=1} \left ( \frac{\sqrt{c}(1+c)}{\xi(\xi+\sqrt{c})}+\frac{c}{\xi+\sqrt{c}}+\frac{c}{\xi^2(\xi+\sqrt{c})}\right )  \, \dif \xi\\
			&=& 2(1+c-(1+c)+c+1-1)\\
			&=& 2c\text{.}
\end{eqnarray*}

\section{Conclusions}
In this paper, we propose a bias-corrected estimator of the noise variance for PPCA models in the high-dimensional framework. The main appeal of our estimator is that it is developed under the assumption that $p/n\to c>0$ as $p, n \to \infty$ and is thus appropriate for a wide range of large-dimensional datasets. Extensive Monte-Carlo experiments demonstrated the superiority of the proposed estimator over several existing estimators (notice however no theoretical justification has been proposed in the literature for these estimators). Moreover, by implementing the proposed estimator of the noise variance within two well-known determination algorithms for the number of principal components, we demonstrate that significant improvement can be obtained. In an additional application and using
this new estimator, we construct the asymptotic theory of the goodness-of-fit test for high-dimensional PPCA models. The overall manage from the paper is that in a high-dimensional PPCA model, when an estimator of the noise variance $\sigma^2$ is needed, the bias-corrected estimator $\widehat{\sigma}_\ast^2$ from the paper should be recommended.

\bibliographystyle{plainnat}
\bibliography{ref,ref2}

\end{document}